%% file: posterior_projector.tex
\documentclass[preprint]{imsart}

\usepackage{amsthm,amsmath,amssymb,graphicx,enumitem,afterpage,bbm}
\usepackage{comment}
\RequirePackage[numbers]{natbib}
\RequirePackage{hyperref}
\hypersetup{colorlinks,
            linkcolor=ultramarine,
            linktoc=ultramarine,
            citecolor=ultramarine,
            urlcolor=ultramarine,
            filecolor=ultramarine
            }

\doi{10.1214/18-EJS1451}

\startlocaldefs
\newcounter{example}[section]
\numberwithin{example}{section}
\newcounter{remark}[section]
\numberwithin{remark}{section}

\newtheorem{theorem}{Theorem}[section]
\newtheorem{corollary}[theorem]{Corollary}
\newtheorem{lemma}[theorem]{Lemma}

\newtheorem{exmp}[example]{Example}
\newtheorem{rmrk}[remark]{Remark}
\newenvironment{example}{\begin{exmp}\rm}{\end{exmp}}
\newenvironment{remark}{\begin{rmrk}\rm}{\end{rmrk}}
\numberwithin{equation}{section}      
\usepackage[math]{easyeqn}

\input mydef
\input ./source/Definitions.tex

\endlocaldefs

\begin{document}

\begin{frontmatter}

\title{Bayesian inference for spectral projectors of the covariance matrix
}
\runtitle{Bayesian inference for spectral projectors}

\begin{aug}
\author{\fnms{Igor} \snm{Silin}
\thanksref{a,b,d,e}
\ead[label=e1]{siliniv@gmail.com}}
%
\and
	\author{\fnms{Vladimir} \snm{Spokoiny}
		\thanksref{f,b,d,e}
		\ead[label=e3]{spokoiny@wias-berlin.de}}	
	
	\address[a]{Moscow Institute of Physics and Technology, \\ 141701, Dolgoprudny, RF.\\
	\printead{e1}}
	\address[f]{Weierstrass Institute and Humboldt University, \\ Mohrenstr. 39, 10117 Berlin, Germany.\\
	\printead{e3}}
	
	\address[b]{National Research University Higher School of Economics,\\ 20 Myasnitskaya ulitsa, 101000, Moscow, RF.}

	\address[d]{Skolkovo Institute of Science and Technology (Skoltech),\\ 143026, Moscow, RF.}	
	
	\address[e]{Institute for Information Transmission Problems RAS,\\ Bolshoy Karetny per. 19, 127051, Moscow, RF.}
		
\end{aug}

\runauthor{I. Silin and V. Spokoiny }

\begin{abstract}
	Let \( X_1, \ldots , X_n \) be an i.i.d. sample in \( \R^p \) with zero mean and the covariance matrix \( \St \). 
	The classical PCA 
	approach recovers the projector \( \Pt \) onto the principal eigenspace of \( \St \) by its empirical counterpart \( \Pe \).
	Recent paper \cite{Koltchinskii_NAACOSPOSC} investigated the asymptotic distribution of the Frobenius distance between the projectors \( \|\Pe - \Pt\|_2 \), while \cite{Naumov_BCSFSPOSC} offered a bootstrap procedure 
to measure uncertainty in recovering 
	this subspace \( \Pt \) even in a finite sample setup.
The present paper considers this problem from a Bayesian perspective and 
suggests to use the credible sets of the pseudo-posterior distribution on the space 
of covariance matrices induced by the conjugated Inverse Wishart prior as sharp confidence sets. 
This yields a numerically efficient procedure.
Moreover, we theoretically justify this method and
derive finite sample bounds on the corresponding coverage probability. 
Contrary to \cite{Koltchinskii_NAACOSPOSC,Naumov_BCSFSPOSC},
the obtained results are valid for non-Gaussian data: the main assumption that we impose is the concentration of the sample covariance \( \Se \) in a vicinity of \( \St \).
	Numerical simulations illustrate good performance of the proposed procedure even on non-Gaussian data in a rather challenging regime.
\end{abstract}





\end{frontmatter}

	\input pp_main_march.tex
	\input ./source/Numerical.tex
	\input pp_proofs_march.tex

	\bibliographystyle{apalike}

\end{document}

%% file: mydef.tex
\renewcommand{\(}{$\,}
\renewcommand{\)}{\,$}

\def\nquad{\hspace{-1cm}}
\def\eqdef{\stackrel{\operatorname{def}}{=}}

\renewcommand{\hat}[1]{\widehat{#1}}
\renewcommand{\tilde}[1]{\widetilde{#1}}

\renewcommand{\Gamma}{\varGamma}
\renewcommand{\Pi}{\varPi}
\renewcommand{\Sigma}{\varSigma}
\renewcommand{\Delta}{\varDelta}
\renewcommand{\Lambda}{\varLambda}
\renewcommand{\Psi}{\varPsi}
\renewcommand{\Phi}{\varPhi}
\renewcommand{\Theta}{\varTheta}
\renewcommand{\Omega}{\varOmega}
\renewcommand{\Xi}{\varXi}
\renewcommand{\Upsilon}{\varUpsilon}

\usepackage{color}

\definecolor{blue(pigment)}{rgb}{0.2, 0.2, 0.6}
\definecolor{ultramarine}{rgb}{0.07, 0.04, 0.56}
\definecolor{darkspringgreen}{rgb}{0.09, 0.45, 0.27}
\definecolor{hookersgreen}{rgb}{0.0, 0.44, 0.0}
\definecolor{plum(traditional)}{rgb}{0.56, 0.27, 0.52}
\definecolor{purple(html/css)}{rgb}{0.5, 0.0, 0.5}
\definecolor{magenta(dye)}{rgb}{0.79, 0.08, 0.48}

%% file: source/Definitions.tex
\newcommand{\data}{\boldsymbol{X}^n}
\newcommand{\subspaceset}{\mathcal{J}}
\newcommand{\Su}{\boldsymbol{\Sigma}}
\newcommand{\Se}{\widehat{\Su}}
\newcommand{\St}{{\Su^{*}}}
\newcommand{\Sp}{\widetilde{\Su}}

\newcommand{\G}{\boldsymbol{G}}

\def\Proj{\boldsymbol{P}}
\def\Projs{\Proj^{*}}
\def\Projh{\widehat{\Proj}}
\def\Projt{\widetilde{\Proj}}

\newcommand{\Pt}{\Projs_{\subspaceset}}
\newcommand{\Ptr}{\Projs_{r}}
\newcommand{\Pts}{\Projs_{s}}
\newcommand{\Pe}{\Projh_{\subspaceset}}

\newcommand{\Per}{\Projh_{r}}

\newcommand{\Pu}{\Proj_{\subspaceset}}
\newcommand{\Pur}{\Proj_{r}}
\newcommand{\Pp}{\Projt_{\subspaceset}}

\def\Gammas{\Gamma^*}
\def\GammasJ{\Gammas_{\subspaceset}}
\def\Gammat{\widetilde{\Gamma}}
\def\GammatJ{\Gammat_{\subspaceset}}
\def\Gammah{\widehat{\Gamma}}
\def\GammahJ{\Gammah_{\subspaceset}}

\def\Lh{\widehat{L}}
\def\LhJ{\Lh_{\subspaceset}}

\def\xiJ{\xi_{\subspaceset}}
\def\xitJ{\widehat{\xi}_{\subspaceset}}
\def\xit{\widehat{\xi}}

\def\Sop{{S}}
\def\Stop{\hat{\Sop}}
\def\StJ{\Stop_{\subspaceset}}

\def\Eu{\boldsymbol{E}}
\def\Ru{\boldsymbol{R}}

\def\IndS{\mathcal{I}_{\subspaceset}}

\newcommand{\Ut}{\boldsymbol{U}^{*}_{\subspaceset}}
\newcommand{\Vt}{\boldsymbol{V}^{*}_{\subspaceset}}
\newcommand{\I}{\boldsymbol{I}}
\def\cond{\, \big| \,}
\def\Cond{\, \bigg| \,}

\def\prior{\Pi}
\newcommand{\Ppost}[1]{\prior\left( #1 \cond \data \right)}

\newcommand{\Pro}{\mathbb{P}}
\newcommand{\R}{\mathbbm{R}}
\newcommand{\SPD}{\mathbbm{S}^{p}_{+}}
\newcommand{\E}{\mathbbm{E}}
\newcommand{\V}{\operatorname{Var}}

\newcommand{\sumn}{\sum\limits_{\textit{j}=1}^\textit{n}}

\newcommand{\ND}{\mathcal{N}}
\newcommand{\rank}{\operatorname{rank}}
\newcommand{\diag}{\operatorname{diag}}
\newcommand{\eff}{\operatorname{eff}}

\def\eqdist{\stackrel{d}{=}}

\def\IWsht{\mathcal{IW}}

\def\erro{\overline{\Diamond}}
\def\erroR{\hat{\Delta}}
\def\erroD{\Delta}
\def\erroS{\mathcal{R}_{\subspaceset}}
\def\hdelta{\widehat{\delta}}

\def\Ind{\operatorname{1}\hspace{-4.3pt}\operatorname{I}}
\def\Var{\operatorname{Var}}
\def\Cov{\operatorname{Cov}}
\def\Fr{2}
\def\T{\top}
\def\Tr{\operatorname{Tr}}

\def\Deltas{\Delta^*}
\def\tdelta{\tilde{\delta}}

\def\ms{m^{*}}
\def\gs{g^{*}}
\def\gh{\hat{g}}
\def\us{u^{*}}
\def\uh{\widehat{u}}

\def\mus{\mu^{*}}
\def\sigmas{\sigma^{*}}
\def\sigmah{\widehat{\sigma}}

%% file: pp_main_march.tex
\section{Introduction} \label{Section: Introduction}
    Let the observed data \( \data = (X_{1}, \ldots, X_{n}) \) be 
    a collection of independent identically distributed 
zero-mean random vectors in \( \R^p \) and let \( X \) be a generic random vector with the same distribution.
Denote by \( \St \) its covariance matrix:
\begin{EQA}
	\St 
	&\eqdef &
	\E \left(X X^{\T}\right).
\end{EQA}
Usually one estimates the true unknown covariance by the sample covariance matrix, given by
\begin{EQA}
    \Se
    & \eqdef &
    \frac{1}{n} \sumn X_{j} X_{j}^{\T}.
\end{EQA}
Quantifying the quality of approximation of \( \St \) by \( \Se \) 
is one of the most classical problems in statistics. 
Surprisingly, a number of deep and strong results in this area 
appeared quite recently. 
The progress is mainly due to Bernstein type results
on the spectral norm \( \| \Se - \St \|_{\infty} \) in the random matrix theory, 
see, for instance, \cite{Koltchinskii_CIAMBFSCO, Rigollet, Tropp, Vershynin_ITTNAAORM, Adamczak}. 
It appears that the quality of approximation is of order \( n^{-1/2} \) while 
the dimensionality \( p \) only enters logarithmically in the error bound.
This allows to apply the results even in the cases of very high data dimension. 

Functionals of the covariance matrix also arise in applications frequently.
For instance, eigenvalues are well-studied in different regimes, see \cite{MP,El Karoui,Johnstone_TW,Wang} and many more references therein.
The Frobenius norm and other \( l_r \)-norms of covariance matrix are of great interest in financial applications; see, e.g. \cite{Fan}.

Much less is known about the quality of estimation of a spectral projector
which is a nonlinear functional of the covariance matrix.
However, such objects arise in dimension reduction methods, manifold learning and spectral methods in community detection, see \cite{Fan_PCA} and references therein for an overview of problems where spectral projectors play crucial role.
Special attention should be focused on the Principal Component Analysis (PCA), probably the most famous dimension reduction method.
Nowadays PCA-based methods are actively used in deep networking architecture \cite{Bengio} and finance \cite{Holtz}, along with other applications.
Over the past decade huge progress was achieved in theoretical guarantees for sparse PCA in high dimensions, see \cite{Johnstone,Birnbaum,Berthet_ODSPCHD,Cai_SPCAORAE,Gao_ROPCSPCA}. 

Suppose we fix some set \( \subspaceset \) of eigenspaces of \( \St \) and consider a direct sum of these eigenspaces and the associated true projector \( \Pt \).
Its empirical counterpart is given by \( \Pe \) computed from the sample covariance \( \Se \).
The recent paper \cite{Reiss} presents new bounds on so called excess risk of PCA defined as 
\( \Tr\left[ \St (\Pt - \Pe) \right] \).

This paper focuses on quantification of uncertainty in recovering the spectral
projector \( \Pt \) from its empirical counterpart \( \Pe \).
More precisely, the random quantity of our interest is the squared Frobenius distance between the true projector and the sample one \( \| \Pe - \Pt \|_{2}^{2} \). 
Even though the projector \( \Pt \) is a complex non-linear mapping of \( \St \), a recent technique from \cite{Koltchinskii_AACBFBFOSPOSC} allows to approximate \( (\Pe-\Pt) \) by a linear functional of \( (\Se-\St) \) with root-n accuracy.
Several results about the distribution of this random variable are available for Gaussian observations: 
\( X_{1}, \ldots, X_{n} \stackrel{i.i.d.}{\sim} \ND(0, \St) \).
For the case when \( \subspaceset \) corresponds to a single eigenvalue in the spectrum and consists of one eigenspace, the normal approximation of \( n\| \Pe - \Pt \|_{2}^{2} \) was shown in \cite{Koltchinskii_NAACOSPOSC} with a tight bound on
\begin{EQA}[c]
    \sup_{x \in \R} 
    \left| \Pro\left\{ 
    \frac{n\| \Pe - \Pt \|_{2}^{2} - \E \left(n\| \Pe - \Pt \|_{2}^{2}\right)}
    	 { \V^{1/2}\left( n\| \Pe - \Pt \|_{2}^{2}\right)} \leq x \right\} 
    - \Phi(x) 
    \right|,
\end{EQA}
where \( \Phi(x) \) is the standard normal distribution function.
However, the distribution of \( n\| \Pe - \Pt \|_{2}^{2} \) depends on the unknown covariance matrix via the first two moments of \( n\| \Pe - \Pt \|_{2}^{2} \).
Due to this fact, it is difficult
to use this result directly for constructing the confidence sets for the true projector \( \Pt \).
The paper \cite{Koltchinskii_NARPCA} demonstrates convergence to a Cauchy-type limit independent from the true covariance operator in a ``high-complexity'' setting where \( \frac{\Tr(\St)}{\| \St \|_{\infty}} \rightarrow \infty \).
Besides, a bootstrap approach can be used to overcome the described problem; see  \cite{Naumov_BCSFSPOSC}.
The bootstrap validity result is based on the  
approximation of the distribution of \( n\| \Pe - \Pt \|_{2}^{2} \)
by the distribution of a Gaussian quadratic form \( \| \xi \|^{2} \).
Namely, for the Gaussian data, Theorem 4.3 of \cite{Naumov_BCSFSPOSC} 
provides the following statement:
\begin{equation}
\begin{aligned}
    \sup_{x\in\R} 
    \left| 
    	\Pro( n\| \Pe - \Pt \|_{2}^{2} \leq x ) - \Pro( \| \xi \|^{2} \leq x)  
    \right| 
    \leq 
    \erro,
\label{Formula: Frequentist approximation}
\end{aligned}
\end{equation}
where \( \xi \) is a zero mean Gaussian vector with a specific covariance 
structure and \( \erro \) is an explicit error term.
The similar approximation is obtained in the bootstrap world, 
this reduces the original problem to the question about Gaussian comparison and Gaussian anti-concentration for large balls.

This paper suggests to look at this problem from a Bayesian point of view.
The standard approach for a nonparametric analysis of the posterior distribution
is based on the prominent Bernstein -- von Mises (BvM) phenomenon.
The BvM result states some pivotal (Gaussian) behavior of the posterior. 
The paper
\cite{Castillo_BvMFSFISM} developed a general framework for functional BvM theorem, while
\cite{Zhou_BvMTFFOCM} used similar ideas to demonstrate asymptotic normality of approximately linear functionals of covariance and precision matrices.
In particular, it can be used to justify the use of Bayesian credible sets 
as frequentist confidence sets for the target parameter; 
see \cite{LeCam,Van der Vaart,Ghosh,Johnstone_BvM,Bickel_BvM,Castillo_NBvMTIGWN} among others. 
In this work, we aim to address a similar question specifically for spectral projectors of the covariance matrix.
It appears that the general BvM technique can be significantly improved and refined for the problem at hand.
The use of the classical conjugated Wishart prior helps not only to build 
a numerically efficient procedure but also to establish precise
finite sample results for the posterior credible sets under mild and general assumptions on the data distribution. 
The key observation here is that,  
similarly to the bootstrap approach of \cite{Naumov_BCSFSPOSC}, 
the credible level sets for the posterior are nearly elliptic,
and the corresponding posterior probability can be approximated by 
a generalized chi-squared-type distribution.
This allows to apply the recent ``large ball probability'' bounds on Gaussian comparison and Gaussian anti-concentration from \cite{Goetze}.  
Moreover, in the contrary to the latter paper \cite{Naumov_BCSFSPOSC}, we do not
require a Gaussian distribution of the data.
Our results claim that the posterior credible sets can be used as frequentist confidence regions even under a possible model misspecification when the true data generation measure is not Gaussian.
We still work with the Gaussian likelihood, in this sense 
our procedure is \emph{pseudo-Bayesian} and the constructed credible sets should be referred to as \emph{pseudo-posterior}; see \cite{Ghosh17}.
In our study we allow the dimension \( p \) to grow with the sample size \( n \), 
however, we need to assume ``\( p^3/n \) is small'' in order to make our results meaningful.

The main contributions of this paper are as follows.

\begin{itemize}
\item 
We offer a new procedure for building elliptic confidence sets for the true projector based on Bayesian simulation from the Inverse Wishart prior.
The procedure is fully data-driven and numerically efficient, 
its complexity is proportional to the squared dimension and independent of sample size.
Numerical simulations confirm good performance of the proposed method for artificial data: both Gaussian and non-Gaussian (not even sub-Gaussian). 

\item
We establish novel results on the coverage properties of {pseudo}-posterior credible sets for a complicated non-linear  
problem of recovering the eigenspace of the sample covariance matrix. 
The results apply under mild conditions on the data distribution.
In particular, we do not require Gaussianity of the observations.
\end{itemize}

The rest of the paper is structured as follows. 
Some notations are introduced in Section~\ref{Subsection: Notations}.
Section~\ref{Subsection: Setup and problem} discusses the model.
Our pseudo-Bayesian framework and the main result of the paper about the pseudo-posterior credible sets are described in 
Section~\ref{Subsection: Bayesian framework and main result}. 
The use of such sets as frequentist confidence sets is discussed in 
Section~\ref{Gaussapprspectr}.
Some numerical results on simulated data are demonstrated in Section~\ref{Section:Numerical experiments}.
Section~\ref{Section: Proofs} contains the proofs of the main theorems.
Some auxiliary results from the literature and the rest of the proofs are collected in Appendix~\ref{Appendix: A} and Appendix~\ref{Appendix: B}, respectively.

\section{Problem and main results} \label{Section: Procedure and main results}
    This section explains our setup and states the main results.

\subsection{Notations} 
\label{Subsection: Notations}

We will use the following notations throughout the paper.
The space of real-valued \( p \times p \) matrices is denoted by \( \R^{p \times p} \), 
while \( \SPD \) means the set of positive-semidefinite matrices. 
We write \( \I_d \) for the identity matrix of size \( d \times d \),
\( \rank(A) \) and \( \Tr(B) \) stand for the \emph{rank} of a matrix \( A \) and the \emph{trace} of a square matrix \( B \).
Further, \( \| A \|_{\infty} \) stands for the \emph{spectral norm} of a matrix \( A \), 
while \( \| A \|_{1} \) means the \emph{nuclear norm}.
The \emph{Frobenius scalar product} of two matrices \( A \) and \( B \) of the same size is \( \langle A, B \rangle_{2} \eqdef \Tr(A^{\T}B) \),  
while the \emph{Frobenius norm} is denoted by \( \| A \|_{2} \). 
When applied to a vector, \( \| \cdot \| \) means just its \emph{Euclidean norm}.
The \emph{effective rank} of a square matrix \( B \) is defined by \( r(B) \eqdef \frac{\Tr(B)}{\| B\|_{\infty}} \).
The relation \( a \lesssim b \) means that there exists an absolute constant \( C \), different from line to line, such that \( a \leq Cb \), while
\( a \asymp b \) means that \( a \lesssim b \) and \( b \lesssim a \).
By \( a\lor b \) and \( a\land b \) we mean maximum and minimum of \( a \) and \( b \), respectively.
In the sequel we will often be considering intersections of events of probability greater than \( 1-{1}/{n} \). 
Without loss of generality, we will write that the probability measure of such an intersection is \( 1 - {1}/{n} \), since it can be easily achieved by adjusting constants.
%
Throughout the paper we assume that \( p < n \).
 
    \subsection{Setup and problem} \label{Subsection: Setup and problem}
Let \( X_{1},\ldots,X_{n} \) be i.i.d. zero mean with \( \Var(X) = \St \).
Without loss of generality, we can assume that \( \St \in \SPD \) is invertible, 
otherwise one can easily transform the data in such a way that the covariance matrix for the transformed data will be invertible. 
Let \( \sigmas_{1} \geq \ldots \geq \sigmas_{p} \) be the ordered eigenvalues of \( \St \).
Suppose that among them there are \( q \) distinct eigenvalues \(  \mus_{1} > \ldots > \mus_q  \).
Introduce groups of indices \( \Deltas_{r} = \{ j: \mus_{r} = \sigmas_{j}\} \) and denote by \( \ms_{r} \) the multiplicity factor (dimension) \( |\Deltas_{r}| \)  for all \( r \in \{1, \ldots, q\} \). 
The corresponding eigenvectors are denoted as \( \us_{1},\ldots, \us_{p} \).
We will use the projector on the \( r \)-th eigenspace of dimension \( \ms_{r} \):
\begin{EQA}
    \Ptr 
    & = &
    \sum_{j\in\Deltas_{r}} \us_{j} {\us_{j}}^{\T}
\end{EQA}
and the eigendecomposition
\begin{EQA}  
    \St 
    & = &
    \sum_{j=1}^p \sigmas_{j} \us_{j} {\us_{j}}^{\T} = \sum_{r=1}^q \mus_{r} \left( \sum_{j\in\Deltas_{r}} \us_{j} {\us_{j}}^{\T} \right) =
        \sum_{r=1}^q \mus_{r} \Ptr.
\end{EQA}
We also introduce the spectral gaps \( \gs_r \):
\begin{EQA}[c]
        \gs_r = 
        \begin{cases}
         \mus_{1} - \mus_{2}, & r = 1,\\
         (\mus_{r-1} - \mus_{r}) \; \land \; (\mus_r - \mus_{r+1}),\;\;& r  \in \{2, \ldots, q-1\},\\
         \mus_{q-1} - \mus_{q}, & r = q.
        \end{cases}     
\end{EQA}  
Suppose that \( \Se \) has \( p \) eigenvalues \( \sigmah_{1} > \ldots > \sigmah_{p} \) (distinct with probability one).
The corresponding eigenvectors are denoted as \( \uh_{1},\ldots, \uh_{p} \).
Suppose that \( \| \Se - \St\|_{\infty} \leq \frac{1}{4} \min\limits_{r\in\{1,\ldots,q\}} \gs_r \).
Then, as shown in \cite{Koltchinskii_AACBFBFOSPOSC}, we can identify clusters of the eigenvalues of \( \Se \) corresponding to each eigenvalue of \( \St \) 
and therefore determine \( \Deltas_r \) and \( \ms_{r} \) for all \( r \in \{1,\ldots, q\} \).
Then we can define the sample projector on the \( r \)-th eigenspace of dimension \( \ms_{r} \):
\begin{EQA}
    \Per 
    & = &
    \sum_{j\in\Deltas_{r}} \uh_{j} {\uh_{j}}^{\T}.
\end{EQA} 
Under the condition that the spectral gap is sufficiently large, 
\cite{Naumov_BCSFSPOSC} approximated
the distribution of \( n\| \Per - \Ptr \|_{2}^{2} \) 
by the distribution of a Gaussian quadratic form \( \| \xi \|^{2} \)
with 
\( \xi~\sim~\ND(0, \Gammas_r) \) and \( \Gammas_r \)
is a block-matrix of the form
\begin{equation}
    \begin{aligned}
        \Gammas_r \eqdef \begin{bmatrix}
        \Gammas_{r1} & O           & \ldots & O         \\
        O           & \Gammas_{r2} & \ldots & O         \\
        \vdots      &  \vdots     & \ddots & \vdots    \\
        O           &  O          & \ldots & \Gammas_{rq}
        \end{bmatrix}
    \label{Formula: Gamma_r_null}
    \end{aligned}
\end{equation}
with \( (q-1) \) diagonal blocks of sizes \( \ms_r \ms_s \times \ms_r \ms_s \):
\begin{EQA}[c]
        \Gammas_{rs} \eqdef \frac{2\mus_{r} \mus_{s}}{(\mus_{r} - \mus_{s})^{2}} \cdot \I_{\ms_r \ms_s}, \;\;\;\; s \neq r.
\end{EQA}
Below we extend these result in two aspects.
First, our approach allows to pick a block of eigenspaces corresponding to 
an interval \( \subspaceset \) in \( \{ 1, \ldots, q \} \) from \( r^{-} \) to \( r^{+} \).
Second, we relax the assumption on Gaussianity of the data.

Let
\begin{EQA}[c]
        \subspaceset = \{ r^-,\; r^- + 1, \; \ldots,\; r^+\}.
\end{EQA}
Define also the subset of indices 
\begin{EQA} 
    \IndS 
    & \eqdef & 
    \bigl\{ k \colon k \in \Deltas_r, \, r \in \subspaceset \bigr\} \,, 
\end{EQA}
and introduce the projector onto the direct sum of the eigenspaces associated with  \( \Ptr \) for all \( r \in \subspaceset \):
\begin{EQA}[c]
    \Pt \eqdef \sum_{r \in \subspaceset} \Ptr 
    = 
    \sum_{k \in \IndS} \us_k {\us_k}^{\T}.
\end{EQA}
Its empirical counterpart is given by
\begin{EQA}[c]
    \Pe \eqdef \sum_{r \in \subspaceset} \Per 
    = 
    \sum_{k \in \IndS} \uh_k \uh_k^{\T}.
\end{EQA}
For instance, when \( \subspaceset = \{ 1, \ldots, q_{\eff}\} \) for some \( q_{\eff} < q \), then \( \Pe \) is exactly what is recovered by PCA. 
Below we focus on \( n\| \Pe - \Pt\|_{2}^{2} \) rather than 
\( n\| \Per - \Ptr \|_{2}^{2} \).

The projector dimension for \( \subspaceset \) is given by \( \ms_{\subspaceset} = \sum_{r \in \subspaceset} \ms_r \). 
Its \emph{spectral gap} can be defined as
\begin{EQA}[c]
        \gs_{\subspaceset} \eqdef \begin{cases}
        \mus_{r^+} - \mus_{r^+ +1}, & \text{ if } r^- = 1;\\
        \mus_{r^- - 1} - \mus_{r^-}, & \text{ if } r^+ = q; \\
        \left( \mus_{r^--1}-  \mus_{r^-}\right) \;\land\; \left( \mus_{r^+} - \mus_{r^+ +1} \right), & \text{ otherwise}. 
        \end{cases}
\end{EQA}
Define also for \( \subspaceset = \{ r^-,\; r^- + 1, \; \ldots,\; r^+\}\)
\begin{EQA}
    l^*_{\subspaceset} 
    & = &
    \mus_{r^-} - \mus_{r^+}.
\end{EQA}
To describe the distribution of the projector \( \Pe \), introduce the following matrix \( \Gammas_{\subspaceset} \) of size \( \ms_{\subspaceset} (p - \ms_{\subspaceset}) \times \ms_{\subspaceset} (p - \ms_{\subspaceset}) \):
\begin{equation} 
    \begin{aligned}
    \Gammas_{\subspaceset} 
    & \eqdef  \diag\left(\Gamma_{\subspaceset}^r\right)_{r \in \subspaceset}, 
    \\
    \Gamma_{\subspaceset}^r 
    & \eqdef 
    \diag\left(\Gamma^{r,s}\right)_{s \notin \subspaceset},
    \\
    \Gamma^{r,s} 
    & \eqdef 
    \frac{2\mus_{r} \mus_{s}}{(\mus_{r} - \mus_{s})^{2}} \cdot \I_{\ms_r \ms_s}, 
    \;\;\;\; r \in \subspaceset, s \notin \subspaceset.
    \label{Formula: Gamma_r}
    \end{aligned}
\end{equation}
It is easy to notice that when \( \subspaceset = \{ r \} \) then this definition coincides with (\ref{Formula: Gamma_r_null}).

Our results apply under one rather mild and natural condition
on the distribution of \( \data = (X_{1}, \ldots, X_{n}) \) that 
the sample covariance matrix \( \Se \) concentrates around the true covariance \( \St \):
\begin{equation}
    \begin{aligned}
        \| \Se - \St\|_{\infty} \leq \hdelta_{n} \|\St\|_{\infty}
\label{Eq: hat_delta_n}
    \end{aligned}
\end{equation}
with probability \( 1 - {1}/{n} \). 
In this result we do not use, say, independence of the \( X_{i} \)'s
or zero mean property, everything is done conditioned on the data \( \data \).
The value \( \hdelta_{n} \) clearly depends on the underlying data distributions, 
but it allows to work with much wider classes of probability measures rather than just Gaussian or sub-Gaussian. 
For the Gaussian case one may take
\begin{EQA}
    \hdelta_{n} 
    & \asymp &  
    \sqrt{\frac{r(\St) + \log(n)}{n}} \,.
\end{EQA}
Several more examples of possible distributions and the corresponding 
\( \hdelta_{n} \) for them are provided in Appendix \ref{Appendix: A}, 
see Theorem \ref{Th: Covariance concentration}. 
So, throughout the rest of the paper we assume that the data satisfy condition (\ref{Eq: hat_delta_n}).

\subsection{Pseudo-Bayesian framework and credible level sets} 
\label{Subsection: Bayesian framework and main result}
Let \( \Pi \) be a prior distribution on the set of
considered covariance matrices \( \Su \).
Even though our data are not necessary Gaussian, 
we can consider the Gaussian log-likelihood:
\begin{EQA}[c]
    \ell_{n}(\Su) 
    = 
    - \frac{n}{2}\log{\det(\Su)} 
    - \frac{n}{2} \Tr(\Su^{-1}\Se) 
    - \frac{np}{2} \log{(2\pi)}.
\end{EQA}
In case of Gaussian data, the posterior measure of a set \( \mathcal{A} \subset \SPD \) can be expressed as
\begin{EQA}[c]
    \Ppost{\mathcal{A}} 
    = 
    \frac{\int_{\mathcal{A}} {\exp{\left(\ell_{n}(\Su) \right)} \;d\Pi(\Su)}}{\int_{\SPD}
    	 {\exp{\left(\ell_{n}(\Su)\right)} \;d\Pi(\Su)}}.
\end{EQA}
However, we can study this random measure for non-Gaussian data as well. 
As the Gaussian log-likelihood \( \ell_{n}(\Su) \) does not necessarily correspond to the true distribution of our data, 
we call the random measure \( \Ppost{\;\cdot\;} \) \emph{a pseudo-posterior}, \cite{Ghosh17}.
Once a prior is fixed, we can easily sample matrices \( \Su \) from this pseudo-posterior distribution.
Denote eigenvalues of \( \Su \) as \( \sigma_{1} > \ldots > \sigma_{p} \) (assume they are distinct with probability one) and eigenvectors as \( u_{1},\ldots, u_{p} \).
The corresponding projector onto the \( r \)-th eigenspace of dimension \( \ms_{r} \) is
\begin{EQA}[c]
    \Pur = \sum_{k\in\Deltas_{r}} u_k {u_k}^{\T}.
\end{EQA}
and the projector on the direct sum of eigenspaces associated with \( \Pur \) for \( r \in \subspaceset \) is
\begin{EQA}[c]
    \Pu = \sum_{r\in\subspaceset} \Pur = \sum_{k \in \IndS} u_k u_k^{\T}.
\end{EQA}
In this work we focus on the conjugate prior to the multivariate Gaussian distribution, that is, the Inverse Wishart distribution \( \IWsht_{p}(\G, p+b-1) \) with \( \G \in \SPD\) , \( 0 < b \lesssim p \).
Its density is given by
\begin{EQA}[c]
     \frac{d\Pi(\Su)}{d\Su} 
     \propto 
     \exp{\left( -\frac{2p+b}{2}\log \det(\Su) - \frac{1}{2} \Tr(\G\Su^{-1})\right)}.
\end{EQA}
Some nice properties of the Inverse Wishart prior distribution allow us to obtain the following result which we will use for uncertainty quantification statements in the next section instead of the Bernstein--von Mises Theorem.

\begin{theorem} 
\label{Theorem: BvM projector conjugate}
    Assume that the distribution of the data \( \data = (X_{1}, \ldots, X_{n}) \) fulfills the sample covariance concentration property (\ref{Eq: hat_delta_n}) with \( \hdelta_n \) satisfying
    \begin{equation}
            \begin{aligned}
            & \hdelta_n \leq \frac{\gs_{\subspaceset}}{4\| \St \|_{\infty}} \;\land\; \frac{r(\St)}{p}.
            \nonumber
            \end{aligned}
        \end{equation}
Consider the prior \( \Pi(\Su) \) given by the Inverse Wishart distribution \\
\( \IWsht_{p}(\G, p+b-1) \).
Let \( \xi \sim \ND(0, \Gammas_{\subspaceset}) \) with 
\( \Gammas_{\subspaceset} \) defined by (\ref{Formula: Gamma_r}).\\
Then with probability \( 1 - \frac{1}{n} \)
\begin{EQA}[c]
    \sup_{x\in\R }\left| 
    	\Ppost{ n\| \Pu - \Pe \|_{2}^{2} \leq x } - \Pro( \| \xi \|^{2} \leq x)  
    \right| 
    \lesssim  
    \Diamond,
\end{EQA}
    where
\begin{equation}
\begin{aligned}
    \Diamond = \Diamond(n,p,\St)
    \eqdef 
    \frac{\Diamond_{1} + \Diamond_{2} + \Diamond_{3}}
         {\| \Gammas_{\subspaceset} \|_2^{1/2} \left( \| \Gammas_{\subspaceset} \|_2^2 - \| \Gammas_{\subspaceset} \|_{\infty}^2 \right)^{1/4}} 
             + \frac{1}{n}\, .
\label{Eq: Diamond}
\end{aligned}
\end{equation}
    The terms \( \Diamond_{1} \) through \( \Diamond_{3} \) can be described as
\begin{EQA}
    \Diamond_{1} 
    & \asymp & 
    \left\{ 
        (\log(n) + p) \left( \,\left( 1 + \frac{l^*_{\subspaceset}}{\gs_{\subspaceset}} \right) \frac{\sqrt{\ms_{\subspaceset}} \| \St \|_{\infty}}{\gs_{\subspaceset}} + \ms_{\subspaceset} \right) \| \St \|_{\infty} + \ms_{\subspaceset}\,\| \G \|_{\infty} 
    \right\} \times  \\
    && \hspace{6.95cm} \times \frac{\ms_{\subspaceset} \, \| \St\|_{\infty}}{{\gs_{\subspaceset}}^{2}}  \sqrt{\frac{\log(n) + p}{n}}
    ,\\
    \Diamond_{2}
    & \asymp & 
    \frac{ \| \St \|_{\infty} \,\left(\ms_{\subspaceset} \, \| \St \|_{\infty}^{2} \land \Tr\left({\St}^{2}\right) \right) }{{\gs_{\subspaceset}}^3} \, p\left( \hdelta_n + \frac{p}{n} \right),\\
    \Diamond_{3} 
    & \asymp &
    \frac{ {(\ms_{\subspaceset})}^{3/2} \| \St \|_{\infty} \, \Tr(\St)}{{\gs_{\subspaceset}}^{2}}
    \sqrt{\frac{\log(n)}{n}}.
\end{EQA}
\end{theorem}
\begin{remark}
The bound (\ref{Eq: Diamond}) can be made more transparent 
if we fix \( \St \) and focus on the dependence on \( p, n, \hdelta_{n} \) and the desired subspace dimension \( \ms_{\subspaceset} \) only (freezing the eigenvalues, the spectral gaps and multiplicities of the eigenvalues):
        \begin{equation}
            \begin{aligned}
                \Diamond \;\asymp\; \sqrt{\frac{(\ms_{\subspaceset})^4\,(p^3 + \log^3(n))}{n}} \;\lor\; \ms_{\subspaceset}\,p\,\hdelta_{n},
            \nonumber
            \end{aligned}
        \end{equation}
or, in the sub-Gaussian case,
    \begin{equation}
            \begin{aligned}
                \Diamond \;\asymp\; \sqrt{\frac{(\ms_{\subspaceset})^4\,(p^3 + \log^3(n))}{n}} .
            \nonumber
            \end{aligned}
        \end{equation}
Moreover, in the case of spiked covariance model we expect $\| \Gammas_{\subspaceset} \|_2$ to behave as $\sqrt{p\,\ms_{\subspaceset}}$ which improves the previous bounds to
\begin{equation}
            \begin{aligned}
                \Diamond \;\asymp\; \sqrt{\frac{(\ms_{\subspaceset})^3\,(p^2 + \log^3(n)/p)}{n}} \;\lor\; \sqrt{\ms_{\subspaceset}\,p}\,\hdelta_{n},
            \nonumber
            \end{aligned}
        \end{equation}
and
    \begin{equation}
            \begin{aligned}
                \Diamond \;\asymp\; \sqrt{\frac{(\ms_{\subspaceset})^3\,(p^2 + \log^3(n)/p)}{n}},
            \nonumber
            \end{aligned}
        \end{equation}
respectively.
\end{remark}

\subsection{Gaussian approximation and frequentist uncertainty quantification
for spectral projectors}
\label{Gaussapprspectr}

For the Gaussian data, Theorem 4.3 of \cite{Naumov_BCSFSPOSC} 
provides the explicit error bound (\ref{Formula: Frequentist approximation})
with the error term \( \erro \) of the following form:
\begin{equation}
    \begin{aligned}
	\sup_{x\in\R} 
	&&
    \left| \Pro( n\| \Per - \Ptr \|_{2}^{2} \leq x ) - \Pro( \| \xi \|^{2} \leq x)  \right| 
    \lesssim \erro,
	\\
    \erro = \erro(n,p,\St) & \eqdef  & 
    \frac{\sqrt{\ms_{r}} \Tr(\Gammas_r)}
    	 {\sqrt{\lambda_{1}(\Gammas_r)\lambda_{2}(\Gammas_r)}} \,\,
    \left( \sqrt{\frac{\log(n)}{n}} + \sqrt{\frac{\log(p)}{n}}\right) 
        \\
        && \qquad 
        + \, \frac{\ms_{r}}{{\gs_{r}}^3} \frac{\Tr^3(\St)}{\sqrt{\lambda_{1}(\Gammas_r)\lambda_{2}(\Gammas_r)}}\sqrt{\frac{\log^3(n)}{n}} \, .
    \label{Eq: overline_Diamond_old}
    \end{aligned}
\end{equation}    
The goal of this section is to extend this result to include the case of 
a generalized spectral cluster and of non-Gaussian data.
Before formulating the result, let us introduce the following auxiliary matrices
    \begin{equation}
        \begin{aligned}
        \Ut & \eqdef 
        \left\{ \frac{{\us_{k}}^{\T}}{\sqrt{\mus_{r}}} \right\}_{\substack{r \in \subspaceset \\ \;\;k\in\Deltas_r }} 
        \in \R^{\ms_{\subspaceset} \times p},
        \\
        \Vt 
        & \eqdef 
        \left\{ \frac{{\us_{l}}^{\T}}{\sqrt{\mus_{s}}} \right\}_{\substack{s \notin \subspaceset \\ \;\;l\in\Deltas_s }} 
        \in \R^{(p- \ms_{\subspaceset})  \times p} \,.
    \nonumber
    \end{aligned}
    \end{equation}
Then the following theorem holds.
\begin{theorem} 
\label{Theorem: clt}
    Assume the distribution of the data \( \data = (X_{1}, \ldots, X_{n}) \) fulfills the sample covariance concentration property (\ref{Eq: hat_delta_n})
    with some \( \hdelta_{n} \).
    Suppose additionally that the projections \( \Pt X \) and \( (\I_{p} - \Pt)X \) are independent and 
    the following third moments are finite: 
    \begin{EQA}
         && \E \| \Ut X \|^3 \leq +\infty, \;\;\; \E \| \Vt X \|^3 \leq +\infty.
    \end{EQA}
    Let \( \xi \sim \ND(0, \Gammas_{\subspaceset}) \) with \( \Gammas_{\subspaceset} \) defined by (\ref{Formula: Gamma_r}) .
    Then
\begin{equation}
\begin{aligned}
     \sup_{x \in \R} 
     \left| 
     	\Pro\left( n\| \Pe - \Pt \|^{2}_{2} \leq x\right) - \Pro(\| \xi\|^{2} \leq x) 
     \right| 
     \lesssim 
     \erro,
\nonumber
\end{aligned}
\end{equation}
    where
\begin{equation}
    \begin{aligned}
\label{Eq: overline_Diamond}
	\erro = \erro(n,p,\St)
	&\eqdef 
	\E \| \Ut X \|^3 \; \E \| \Vt X \|^3 \, \frac{p^{1/4}}{\sqrt{n}} \\
	& \quad
	+ \frac{\overline{\Delta}}{\| \Gammas_{\subspaceset} \|_2^{1/2} \left( \| \Gammas_{\subspaceset} \|_2^2 - \| \Gammas_{\subspaceset} \|_{\infty}^2 \right)^{1/4}}\,,
	\qquad
	\\
	\overline{\Delta} 
    & \eqdef 
    n \ms_{\subspaceset} \left( 1 + \frac{l^*_{\subspaceset}}{\gs_{\subspaceset}}\right)\,
    \left\{ \left( 1 + \frac{l^*_{\subspaceset}}{\gs_{\subspaceset}}\right)\frac{\hdelta_{n}^4}{{\gs_{\subspaceset}}^4} \;\lor\; 
            \frac{|\subspaceset| \;\hdelta_{n}^3}{{\gs_{\subspaceset}}^3}
    \right\} \,.
    \end{aligned}
\end{equation}
\end{theorem}
\begin{remark}
    The condition on independence of \( \Pt X \) and \( (\I_{p} - \Pt)X \) is not very restrictive, in fact, it has a natural interpretation:
    while we are interested in the ``signal'' \( \Pt X \), the orthogonal part \( (\I_{p} - \Pt)X \) can be considered as ``noise'', and it is plausible to assume that the ``noise'' is independent from the ``signal''.
    It also worth mentioning that this condition was not required in our main result above about the behavior of the pseudo-posterior.
\end{remark}
\begin{remark}
   The components of \( \Ut X \) are the $\ms_{\subspaceset}$-dimensional coordinates of $X$ after projecting onto the eigenspace of interest and proper scaling.
    Similarly, the components of \( \Vt X \) are the $(p-\ms_{\subspaceset})$-dimensional coordinates of $X$ after projecting onto the orthogonal complement and proper scaling.
    In general, the factors \( \E \| \Ut X \|^3 \) and \( \E \| \Vt X \|^3 \) from the error bound \( \erro \) depend on how heavy the tails of the distribution of \( X \) are.
    However, it is easy to show that in case of a sub-Gaussian random vector \( X \) the behaviour is as follows:
    \begin{EQA}
         \E \| \Ut X \|^3 & \lesssim & (\ms_{\subspaceset})^{3/2}\,,\\
         \E \| \Vt X \|^3 & \lesssim & {(p-\ms_{\subspaceset})}^{3/2}\,.
    \end{EQA}
    Coupled with Theorem \ref{Th: Covariance concentration}, (ii), this allows to bound the error term \( \erro \) in the sub-Gaussian case as
    \begin{EQA}
        \erro & \lesssim & \sqrt{\frac{(\ms_{\subspaceset})^{3} {p}^{3.5}}{n}} + {\ms_{\subspaceset}} |\subspaceset| \sqrt{\frac{p^3 + \log^3(n)}{n}} \, ,
    \end{EQA}
    where, for simplicity, the characteristics of \( \St \) are hidden in the constants. 
\end{remark}
The proof of this result is presented in Appendix \ref{Appendix: B}.
The obtained bound is worse than (\ref{Eq: overline_Diamond_old}) when 
the the full dimension \( p \) is large.
This is the payment for the Gaussian approximation which appears for non-Gaussian data. Our result makes use of the Gaussian approximation technique from \cite{Bentkus_LTB}.
Some recent developments in Gaussian approximation for a probability of a ball indicate that the bound 
(\ref{Eq: overline_Diamond}) can be improved even further;
see \cite{Spdraft2018}.
Comparison of the results of Theorem \ref{Theorem: BvM projector conjugate} and Theorem \ref{Theorem: clt} reveals 
that the pseudo-posterior distribution of  \( n\| \Pu - \Pe \|_{2}^{2} \) given the data
perfectly mimics the distribution of \\ \(  n\| \Pe - \Pt \|_{2}^{2} \),
and, therefore, can be applied to building of elliptic confidence sets for the true projector.
Specifically, for any significance level \( \alpha \in (0;\;1) \) (or confidence level \( 1-\alpha \)) we can estimate the true quantile 
\begin{EQA}[c]
    \gamma_{\alpha} 
    \eqdef 
    \inf \left\{ 
    	\gamma > 0: \Pro\left(n\| \Pe - \Pt\|_{2}^{2} > \gamma \right) \leq \alpha 
    \right\}
\end{EQA}
by the following counterpart which can be numerically assessed using Bayesian 
credible sets:
\begin{EQA}[c]
     \gamma_{\alpha}^{\circ} 
     \eqdef 
     \inf \left\{ 
     	\gamma > 0: \Ppost{n\| \Pu - \Pe\|_{2}^{2} > \gamma\;} \leq \alpha 
     \right\}.
\end{EQA}
Then, the main results presented above imply
the following corollary.

\begin{corollary} \label{Corollary: Confidence set Wishart}
Assume that all conditions of Theorem \ref{Theorem: BvM projector conjugate} and Theorem \ref{Theorem: clt} are fulfilled.
Then
\begin{EQA}[c]
	\sup_{\alpha\in(0;\;1)} 
	\left| 
		\alpha - \Pro\left(n\| \Pe - \Pt\|_{2}^{2} > \gamma_{\alpha}^{\circ}\right) 
	\right| 
	\lesssim 
	\Diamond + \erro,
\end{EQA}
    where \( \Diamond = \Diamond(n,p,\St) \), \( \erro = \erro(n,p,\St) \) are defined by (\ref{Eq: Diamond}), (\ref{Eq: overline_Diamond}), respectively.
\end{corollary}

%% file: source/Numerical.tex
\section{Numerical experiments} \label{Section:Numerical experiments}
	This section shows by mean of artificial data that the proposed pseudo-Bayesian approach works quite well even 
for large data dimension and limited sample size.
We also want to track how the quality depends on the sample size \( n \) and the dimension \( p \).
In our experiments we first  
fix some true covariance matrix \( \St \) of size \( p \times p \).
Without loss of generality we consider only diagonal matrices \( \St \), 
so \( \St \) is defined by the distinct eigenvalues \( \mus_{r} \) and the multiplicities \( \ms_{r} \). 
We also specify the desired subspace that we want to investigate by fixing \( \subspaceset \).
Further, for different sample sizes \( n \) we repeat the following two-step procedure.
The first step is to determine the quantiles of \( n\| \Pe - \Pt \|_{2}^{2} \). 
For that we generate \( 3000 \) samples \( \data \), compute the corresponding \( \Pe \) 
and then just take \( \alpha- \)quantiles of the obtained realizations \( n\| \Pe - \Pt \|_{2}^{2} \) for \( \alpha \) from \( 0.001 \) to \( 0.999 \) with step \( 0.001 \).
The second step is to estimate the quantiles of the pseudo-posterior distribution of \( n\| \Pu - \Pe \|_{2}^{2} \).
We generate \( 50 \) samples \( \data \) and for each realization we generate \( 3000 \) pseudo-posterior covariance matrices \( \Su \) from the Inverse Wishart distribution with \( G= \I_p,\; b = 1 \).
Then we compute the corresponding \( \Pu \) and take the \( \alpha- \)quantiles of \( n\| \Pu - \Pe \|_{2}^{2} \) just as in the first step.
Namely, for each \( \alpha \) we get \( 50 \) quantile estimates \( {\gamma_{\alpha}^{\circ}}^{(j)},\; j \in \{ 1, \ldots, 50\} \) (suppose we order them in ascending order) and take median of them. 
For the true quantiles from the first step and the medians of the quantile estimates from the second step we build a QQ-plot, which consists of points with coordinates \(\left(\gamma_{\alpha}, {\gamma_{\alpha}^{\circ}}^{(25)}\right)\) for various \( \alpha \).
We expect that the constructed QQ-plot is close to the identity line indicating that these two distributions are close to each other. 
Also we present a table with median coverage probabilities 
\[
	\Pro\left(n\| \Pe - \Pt \|_{2}^{2} \leq {\gamma_{\alpha}^{\circ}}^{(25)}\right) 
\] 
and interquartile ranges 
\[
	\Pro\left(n\| \Pe - \Pt \|_{2}^{2} \leq {\gamma_{\alpha}^{\circ}}^{(38)}\right) - \Pro\left(n\| \Pe - \Pt \|_{2}^{2} \leq {\gamma_{\alpha}^{\circ}}^{(12)}\right) 
\]
 of this coverage probability for the desired confidence levels \( 1-\alpha \) from the list
\[
 	 \{0.99, 0.95, 0.90, 0.85, 0.80, 0.75 \}.
\]

\begin{figure}[!htp]
	\center{
		\includegraphics[width=0.45\linewidth]{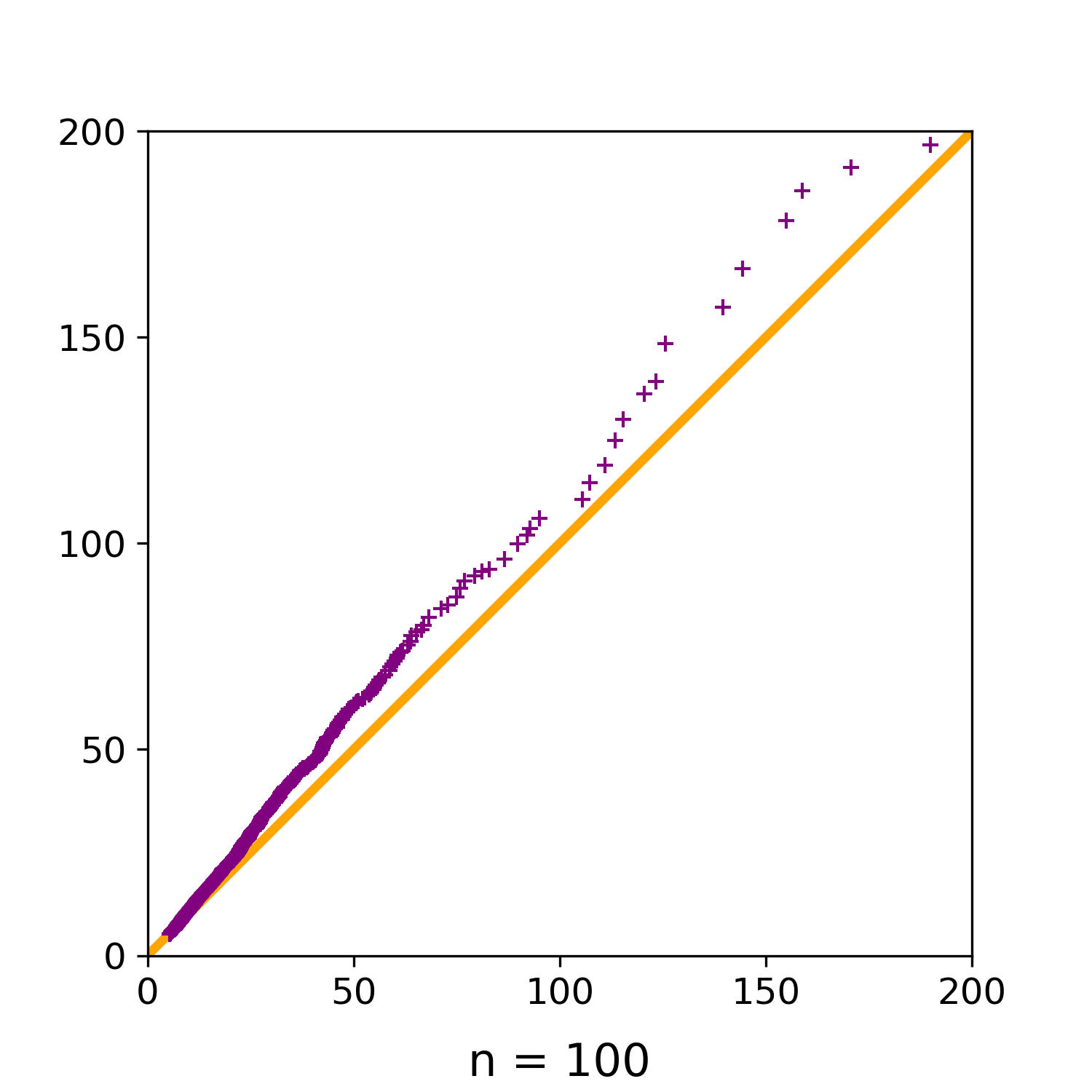}
	\hfill
		\includegraphics[width=0.45\linewidth]{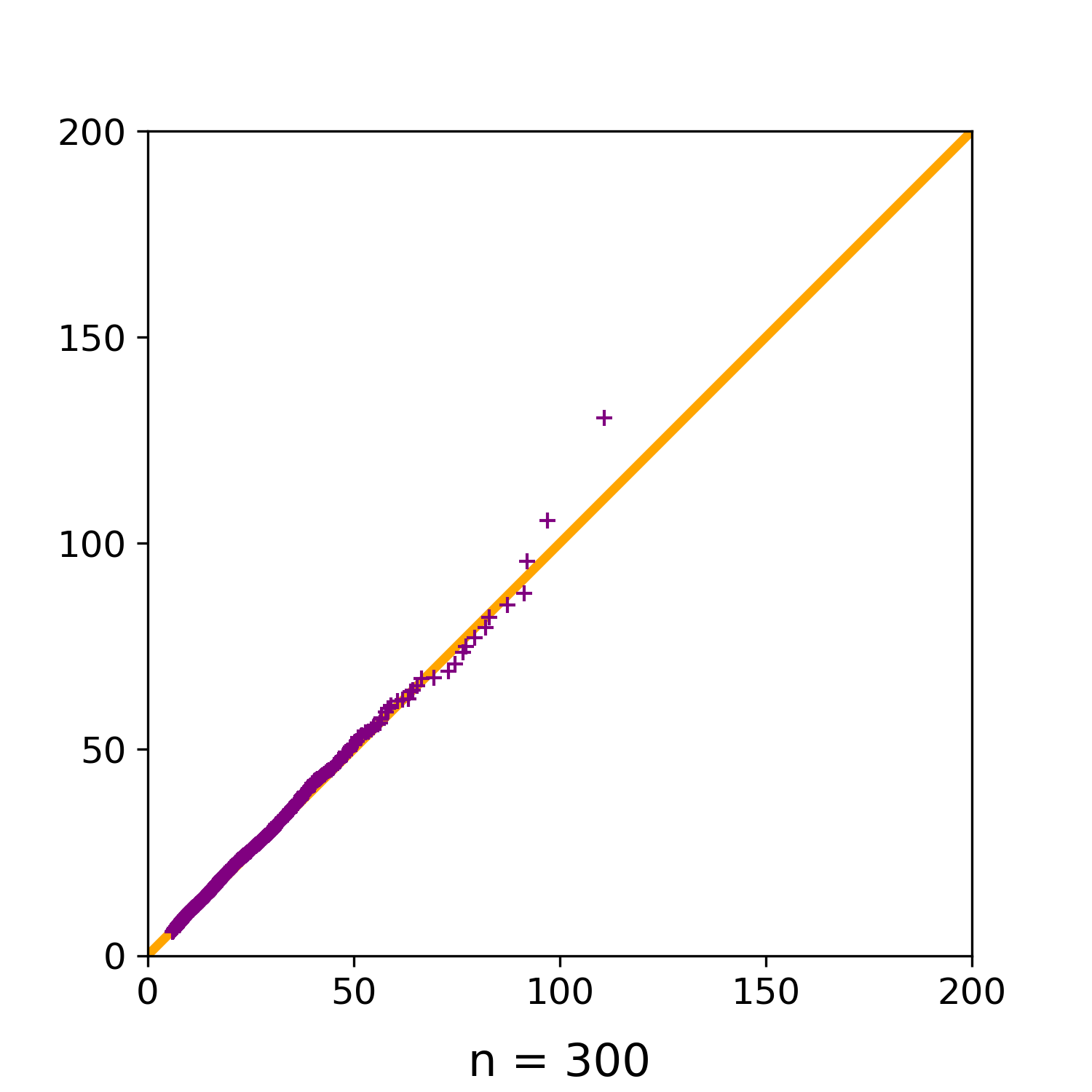}
	} \\
	\center{
		\includegraphics[width=0.45\linewidth]{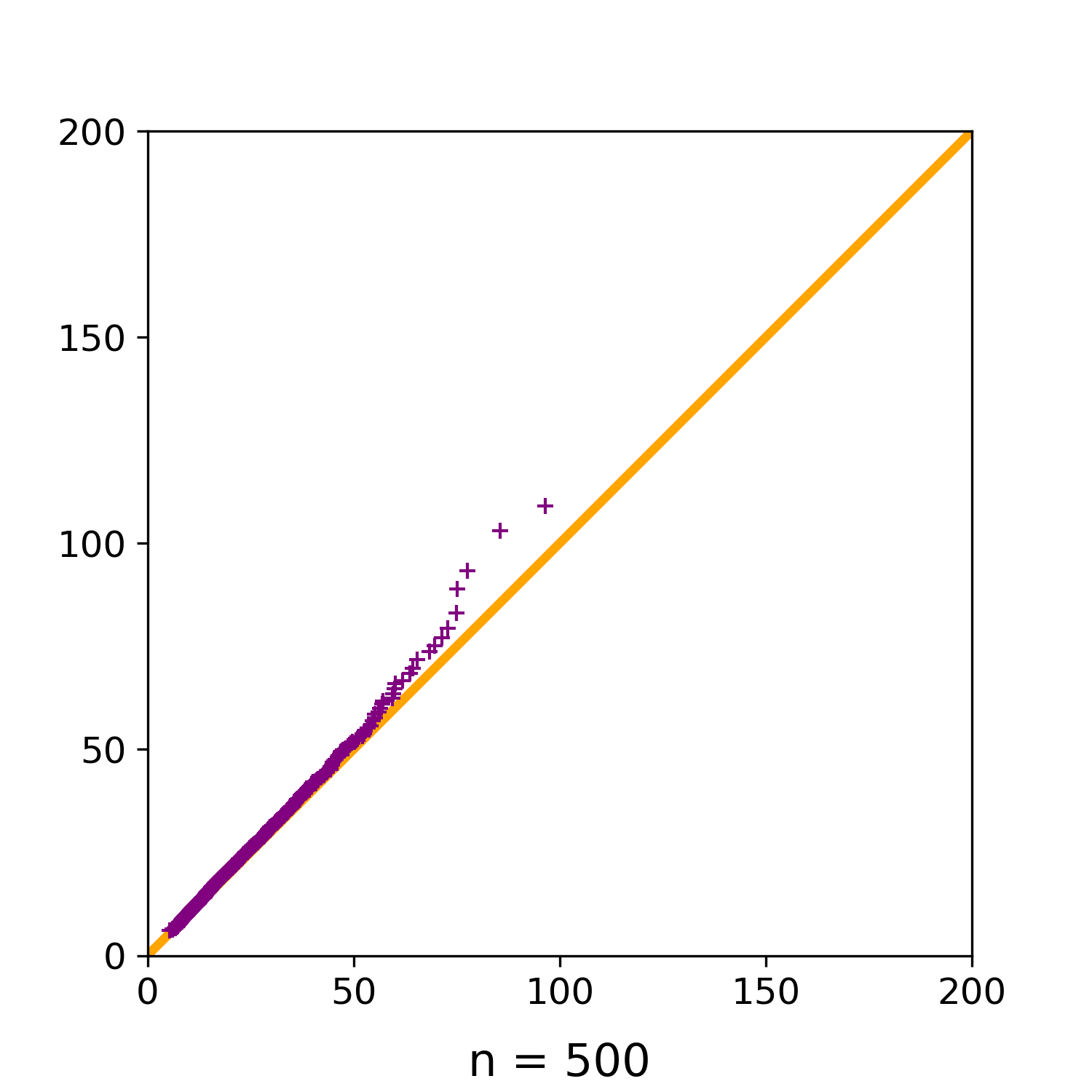}
	\hfill
		\includegraphics[width=0.45\linewidth]{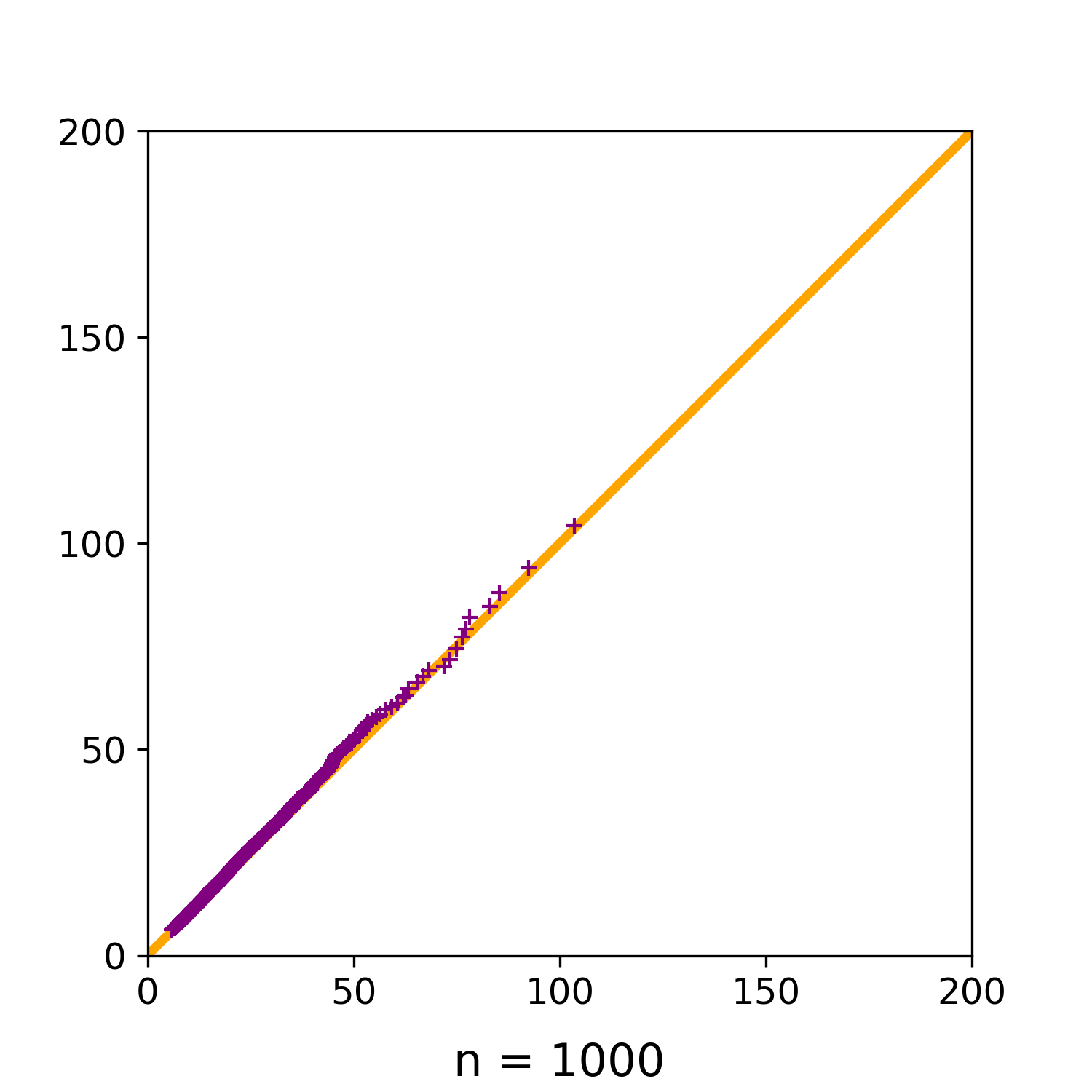}
	} \\
	\center{
		\includegraphics[width=0.45\linewidth]{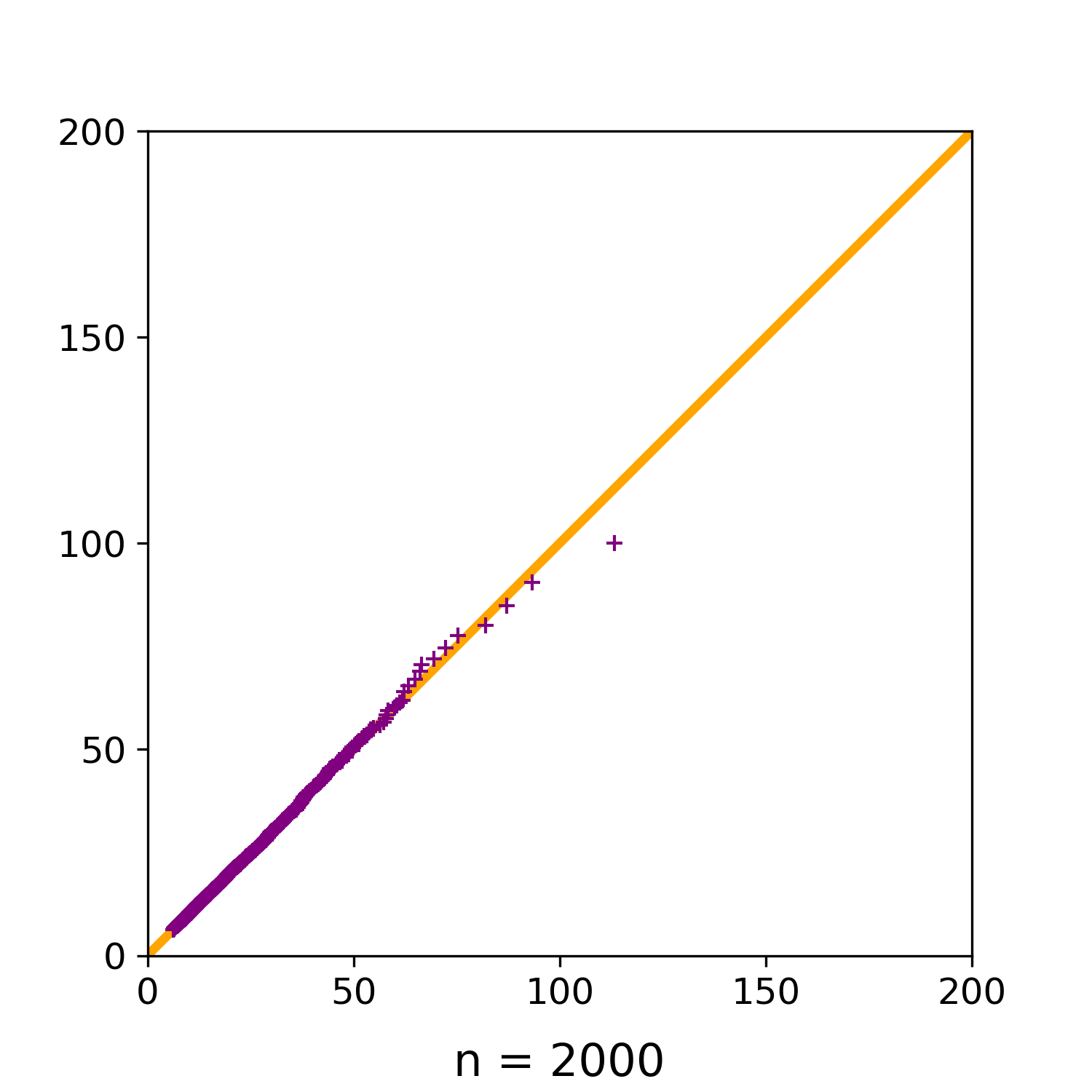}
	\hfill
		\includegraphics[width=0.45\linewidth]{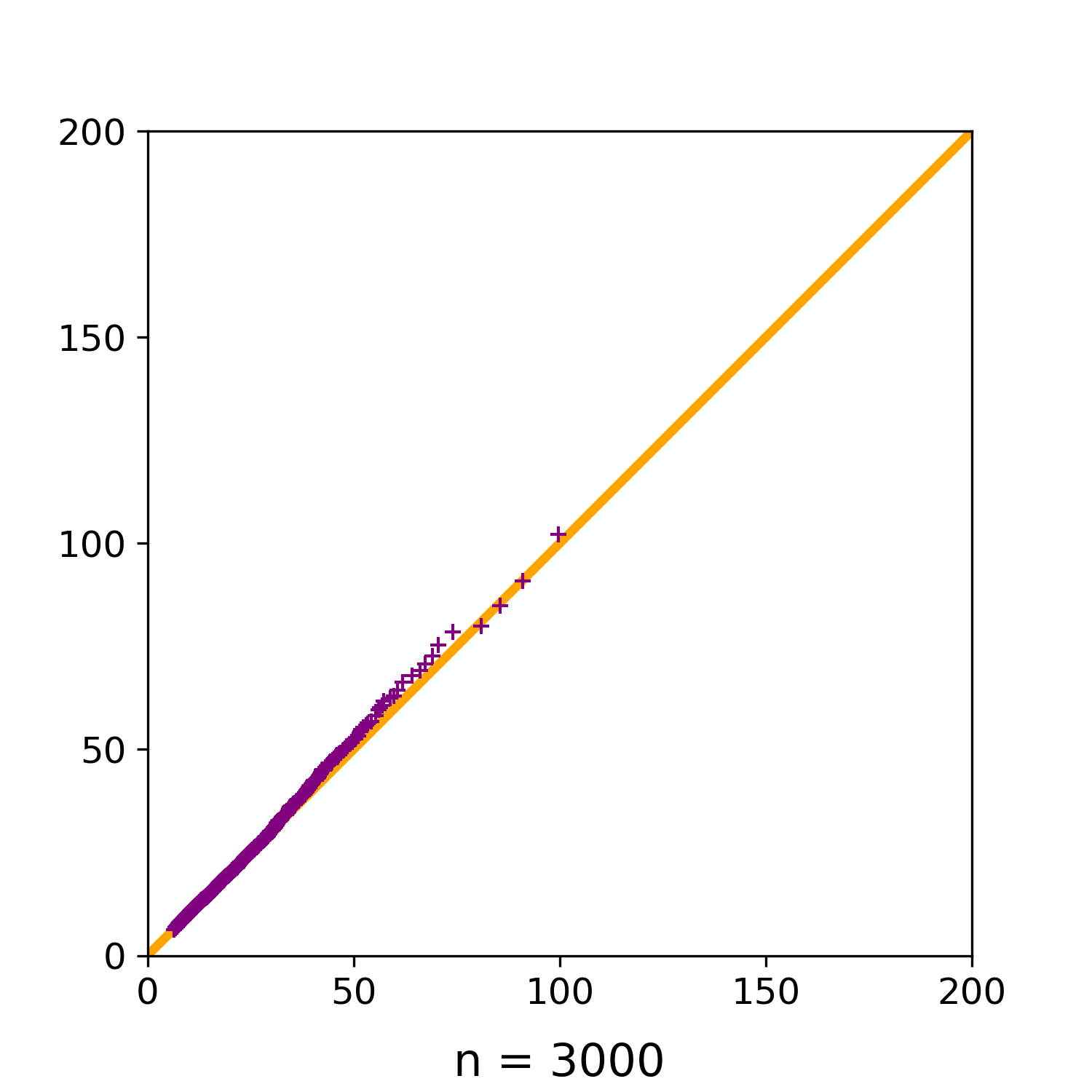}
	}
	\caption[]{QQ-plots of the proposed pseudo-Bayesian procedure for the first experiment (Gaussian data).}
	\label{Fig:1}
\end{figure}

In the first experiment we work with Gaussian data.
The parameters of the experiment are as follows:
\begin{itemize}
	\item \( p=100 \), \( \ms_{r} = 1 \) for all \( r \in \{1,\ldots, 100\} \).
	\item \( \mus_1 = 25.698,\; \mus_{2} = 15.7688,\; \mus_3 = 10.0907,\; \mus_4 = 5.9214, \;\mus_5=3.4321 \) and the rest of the eigenvalues \( \mus_6, \ldots, \mus_{100} \) are from the Marchenko -- Pastur law with support \( [0.71;\;1.34] \).
	\item \( \subspaceset = \{ 1 \} \), so we investigate the one-dimensional principal subspace given by \( {\Projs_{1}} \).
\end{itemize} 
The QQ-plots are depicted on Figure~\ref{Fig:1} while the coverage probabilities and the interquartile ranges are presented in Table~\ref{Tab:1}.

The setup of this experiment is exactly the same as the second example of \cite{Naumov_BCSFSPOSC}, so the performance of our pseudo-Bayesian method can be directly compared with the performance of the Bootstrap approach (cf. Figure 2 and Table 2 of \cite{Naumov_BCSFSPOSC}). The accuracy of these two procedures is approximately the same. 

\begin{figure}[!htp]
	\center{
		\includegraphics[width=0.45\linewidth]{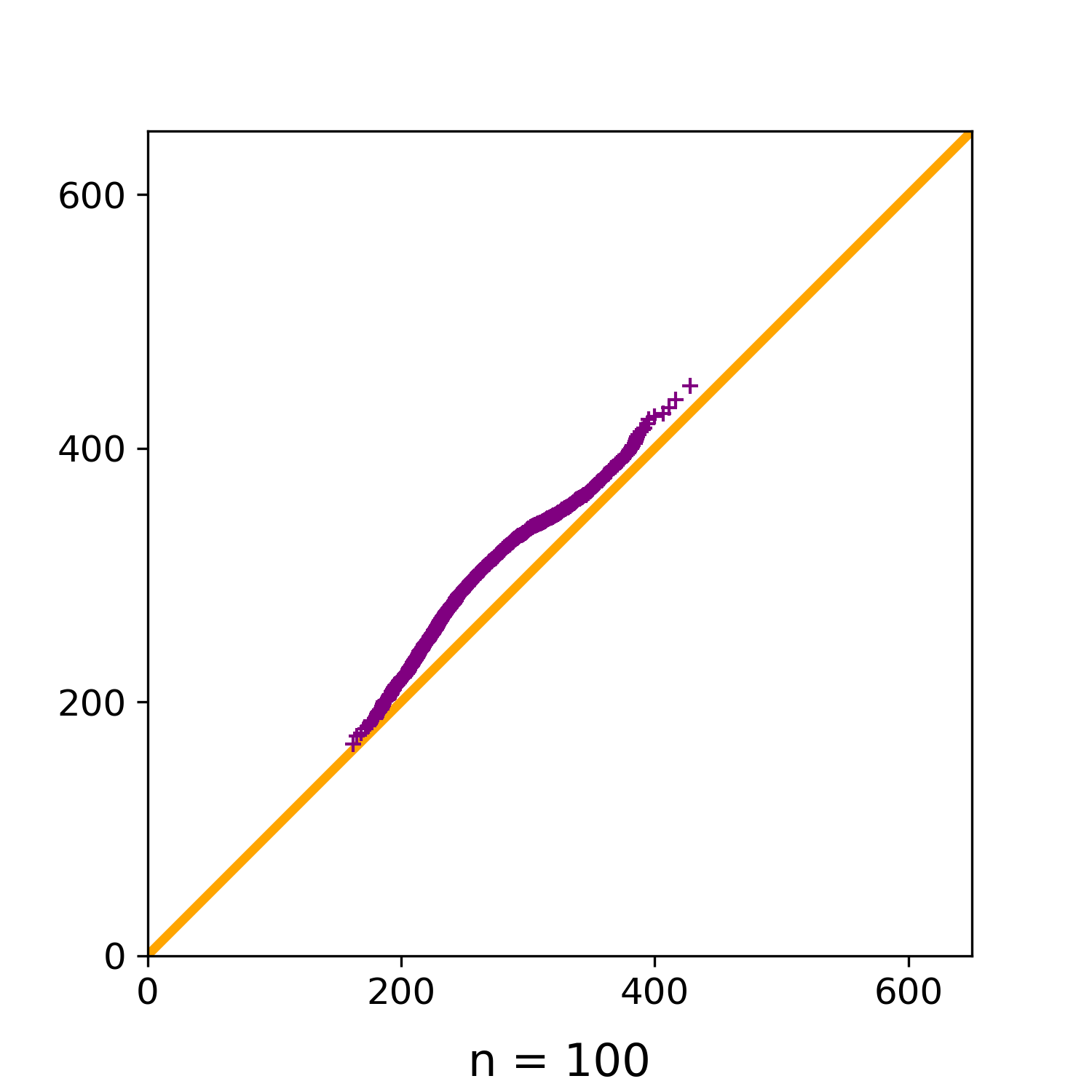}
	\hfill
		\includegraphics[width=0.45\linewidth]{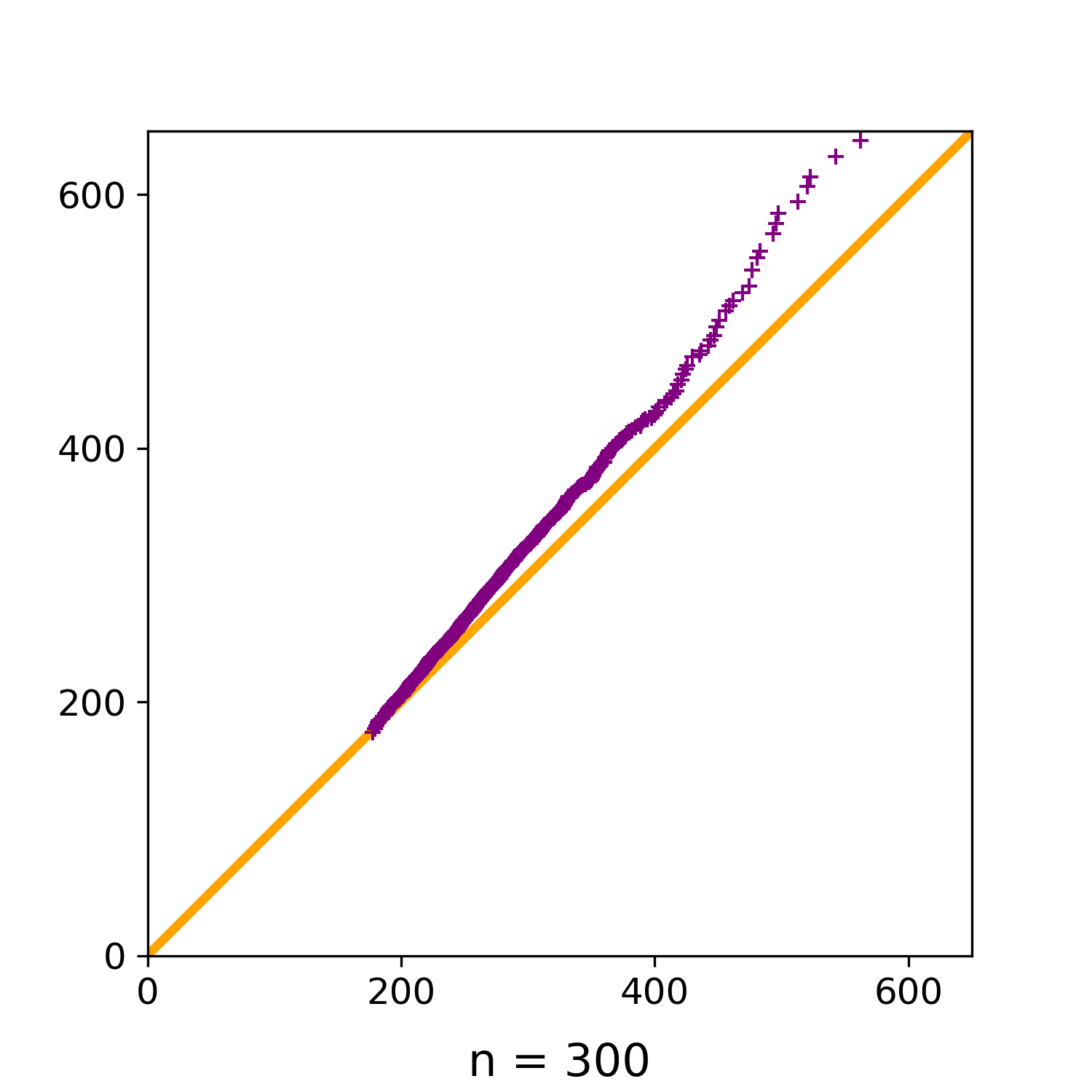}
	} \\
	\center{
		\includegraphics[width=0.45\linewidth]{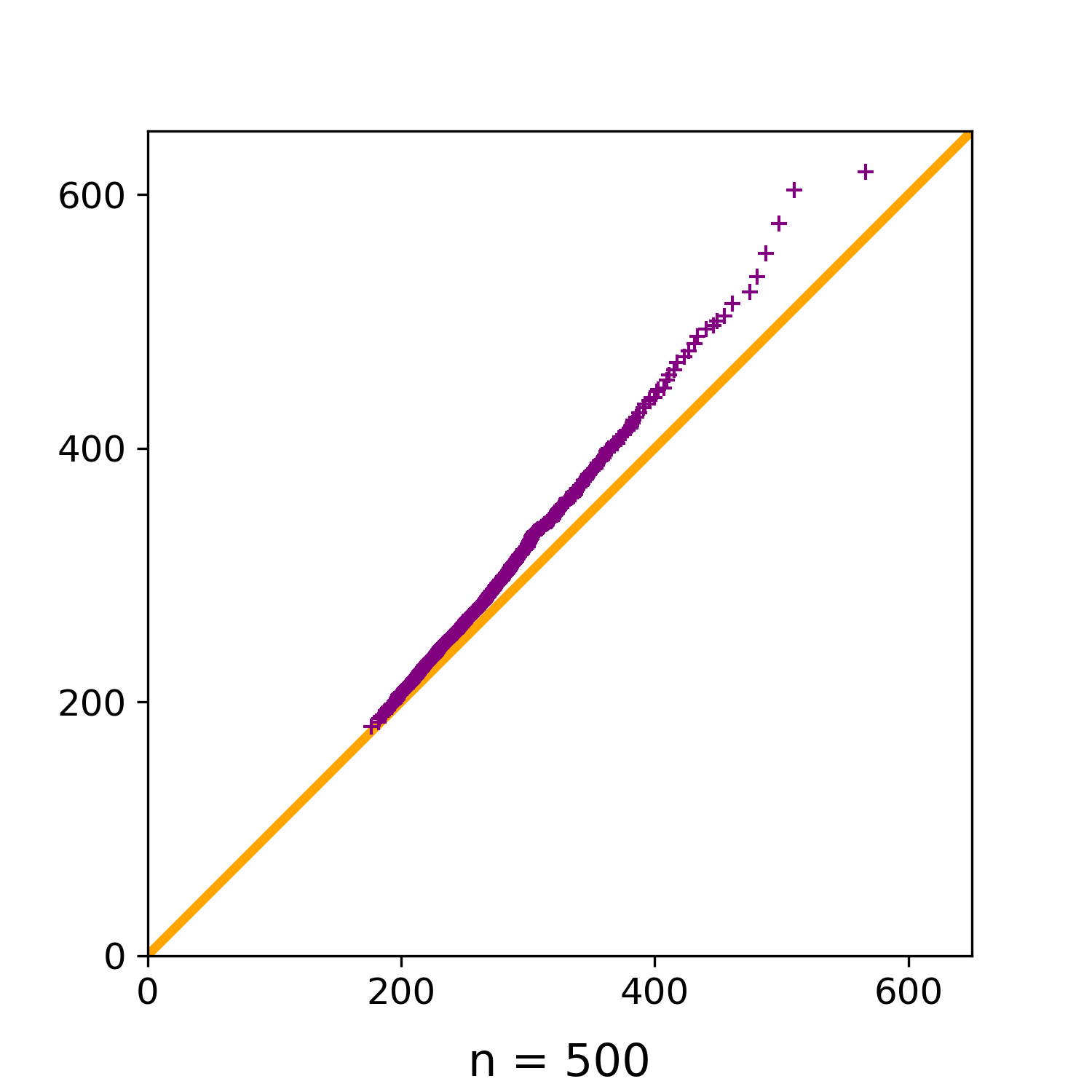}
	\hfill
		\includegraphics[width=0.45\linewidth]{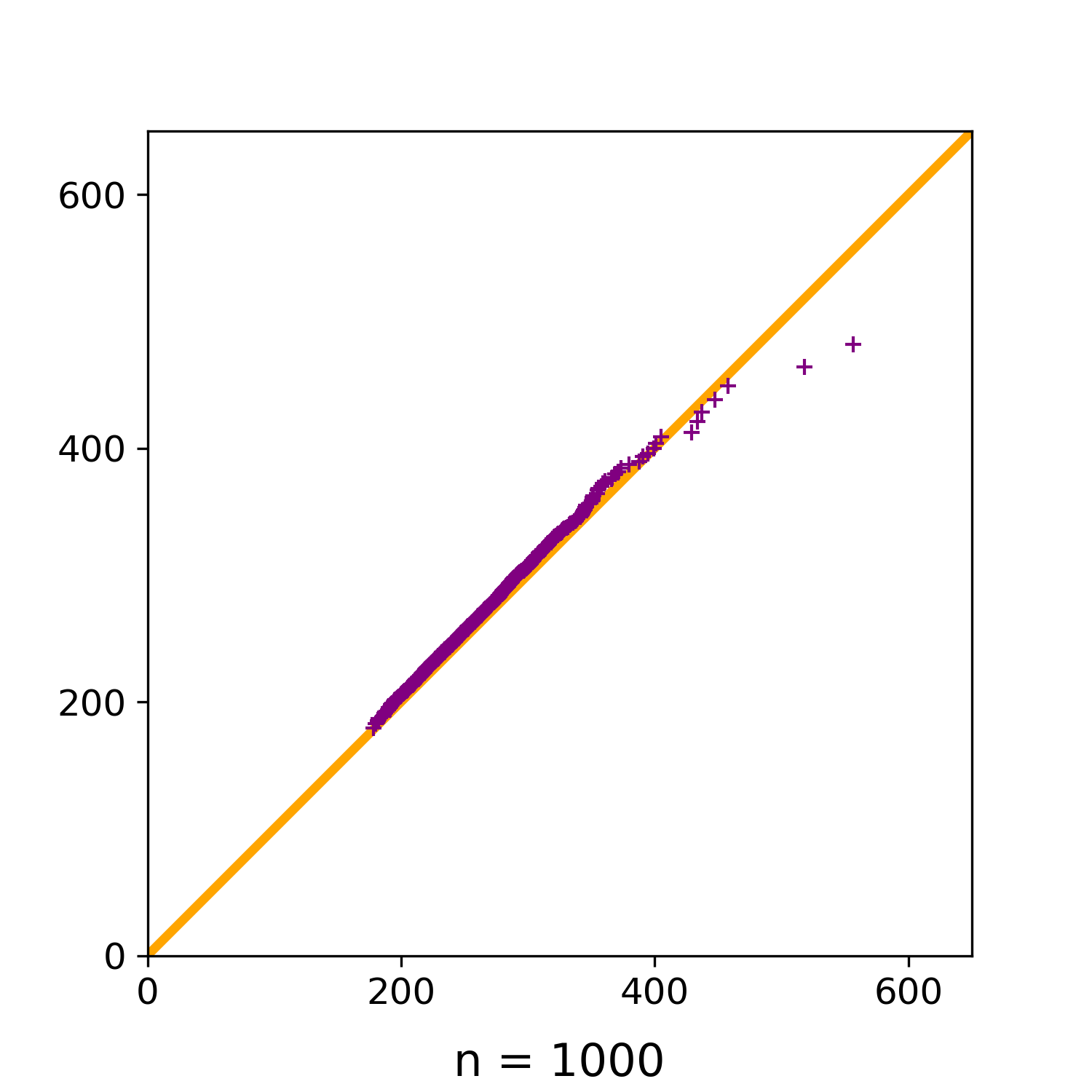}
	} \\
	\center{
		\includegraphics[width=0.45\linewidth]{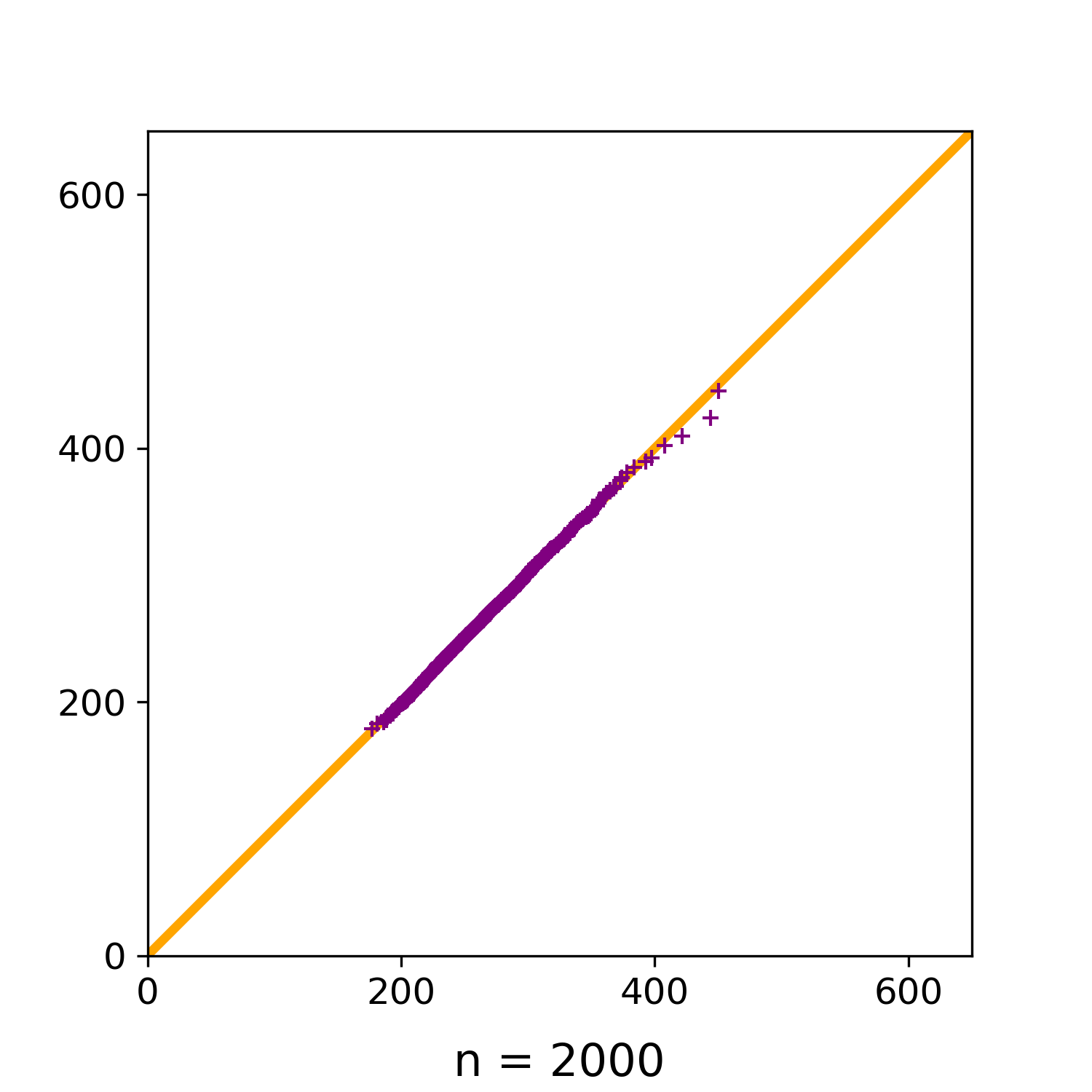}
	\hfill
		\includegraphics[width=0.45\linewidth]{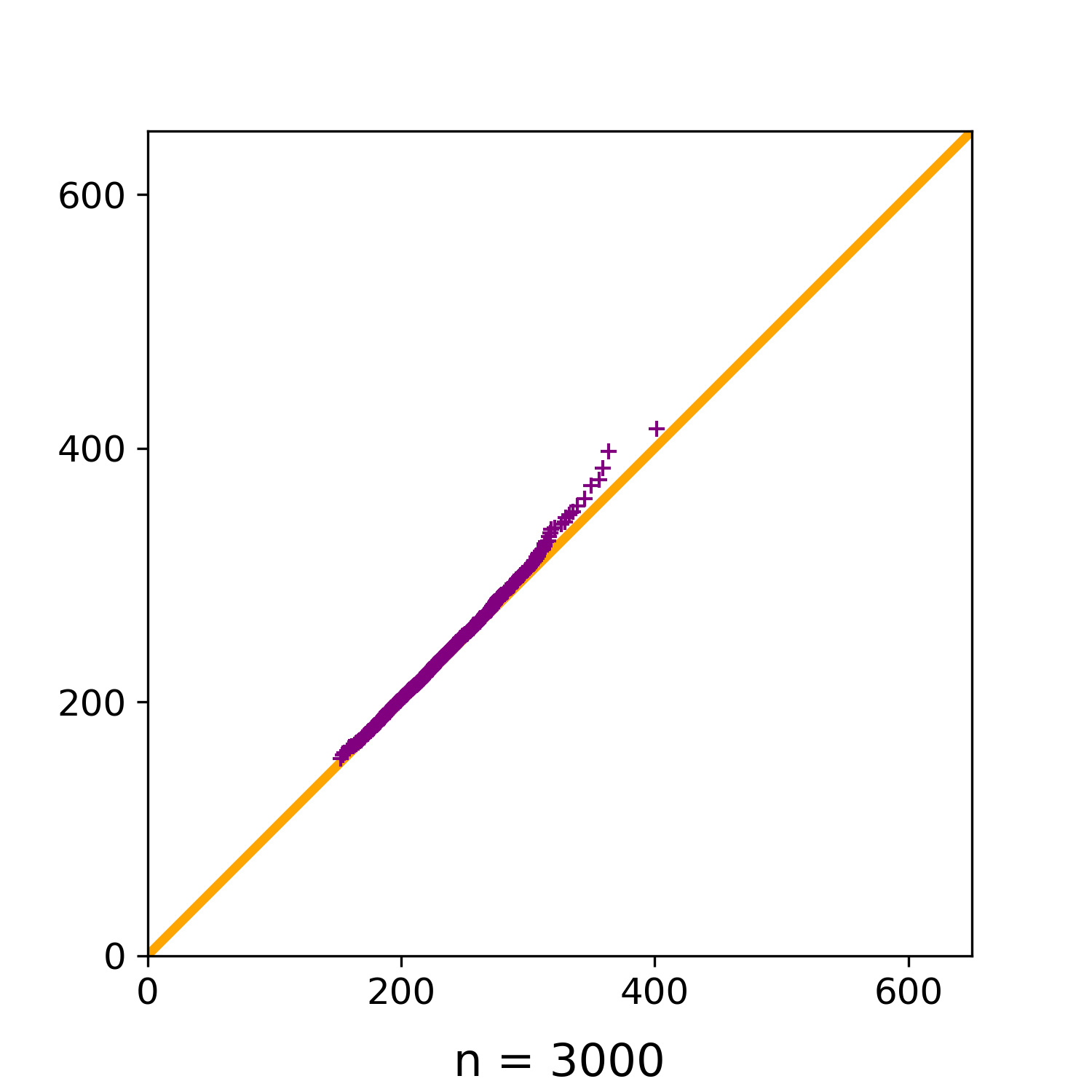}
	}
	\caption{QQ-plots of the proposed pseudo-Bayesian procedure for the second experiment (non-Gaussian data).}
	\label{Fig:2}
\end{figure}

In the second experiment we check how our method performs on non-Gaussian data.
We generate each component of the vectors \( X_j \) independently yielding 
diagonal covariance matrix.
In addition to Gaussian distribution, we consider also the following three options: the uniform distribution on the interval \( [-a;\;a] \), the Laplace distribution with scaling parameter \( a \) and the discrete uniform distribution with three values \( \{ -a, 0, a\} \).
In each case the parameter \( a \) is chosen in such a way that ensures the variance located on the diagonal of the covariance matrix fixed earlier. 
So, the parameters of the experiment are as follows:
\begin{itemize}
	\item \( p = 100 \), \( \ms_1 = 3 \), \( \ms_{2} = 3 \), \( \ms_3 = 3 \) and the rest of the multiplicities \( \ms_4,\ldots, \ms_{91} \) are one.
	\item \( \mus_1 = 25 \), \( \mus_{2} = 20 \), \( \mus_3 = 15 \), \( \mus_4=10 \), \( \mus_5=7.5 \), \( \mus_6=5 \) and the rest of the eigenvalues 
	\( \mus_7, \ldots, \mus_{100} \) are from the uniform distribution on \( [0;\;3] \). 
	\item The first nine components were generated according to: uniform, Laplace, discrete, Gaussian, Laplace, discrete, Laplace, Laplace, uniform distributions, respectively. The rest of the components are Gaussian. 
	\item \( \subspaceset = \{ 1, 2, 3 \} \), so we investigate nine-dimensional subspace given by \\\( \Projs_{1} + \Projs_{2} + \Projs_{3} \).
\end{itemize} 
The QQ-plots are depicted on Figure~\ref{Fig:2} while the coverage probabilities and the interquartile ranges are presented in Table~\ref{Tab:2}.

\begin{table*}[!htp]
	\caption{Coverage probabilities and interquartile ranges of the proposed pseudo-Bayesian procedure for the first experiment (Gaussian data). }
	\label{Tab:1}
	\centering
	\begin{tabular}{ccccccc}
	\hline
	     & \multicolumn{6}{c}{Confidence levels (\( 1-\alpha \))} \\ \cline{2-7} 
	\centering \( n \)    & \( 0.99 \)   & \( 0.95 \)   & \( 0.90 \)   & \( 0.85 \)   & \( 0.80 \)   & \( 0.75 \)  \\ \hline
	\( 100 \)  & \( 0.993 \) & \( 0.968 \) & \( 0.929 \) & \( 0.893 \) & \( 0.854 \) & \( 0.809 \) \\
	       & \( 0.023 \) & \( 0.061 \) & \( 0.103 \) & \( 0.145 \) & \( 0.176 \) & \( 0.204 \) \\ \cline{2-7} 
	\( 300 \)  & \( 0.988 \) & \( 0.952 \) & \( 0.906 \) & \( 0.851 \) & \( 0.805 \) & \( 0.762 \) \\
	       & \( 0.026 \) & \( 0.085 \) & \( 0.143 \) & \( 0.184 \) & \( 0.199 \) & \( 0.216 \) \\ \cline{2-7} 
	\( 500 \) & \( 0.993 \) & \( 0.955 \) & \( 0.909 \) & \( 0.865 \) & \( 0.812 \) & \( 0.771 \) \\
	       & \( 0.022 \) & \( 0.072 \) & \( 0.099 \) & \( 0.108 \) & \( 0.123 \) & \( 0.126 \) \\ \cline{2-7} 
	\( 1000 \) & \( 0.990 \) & \( 0.956 \) & \( 0.908 \) & \( 0.859 \) & \( 0.817 \) & \( 0.767 \) \\
	       & \( 0.014 \) & \( 0.049 \) & \( 0.066 \) & \( 0.067 \) & \( 0.091 \) & \( 0.104 \) \\ \cline{2-7} 
	\( 2000 \) & \( 0.992 \) & \( 0.952 \) & \( 0.898 \) & \( 0.847 \) & \( 0.793 \) & \( 0.747 \) \\
	       & \( 0.009 \) & \( 0.030 \) & \( 0.058 \) & \( 0.064 \) & \( 0.067 \) & \( 0.065 \) \\ \cline{2-7} 
	\( 3000 \) & \( 0.992 \) & \( 0.959 \) & \( 0.908 \) & \( 0.849 \) & \( 0.802 \) & \( 0.750 \) \\
	       & \( 0.005 \) & \( 0.024 \) & \( 0.036 \) & \( 0.054 \) & \( 0.046 \) & \( 0.063 \) \\ \hline
	\end{tabular}
\end{table*}

\begin{table*}[!htp]
	\caption{Coverage probabilities and interquartile ranges of the proposed pseudo-Bayesian procedure for the second experiment (non-Gaussian data) }
	\label{Tab:2}
	\centering
	\begin{tabular}{ccccccc}
	\hline
	     & \multicolumn{6}{c}{Confidence levels (\( 1-\alpha \))} \\ \cline{2-7} 
	\centering \( n \)    & \( 0.99 \)   & \( 0.95 \)   & \( 0.90 \)   & \( 0.85 \)   & \( 0.80 \)   & \( 0.75 \)  \\ \hline
	\( 100 \)  & \( 0.997 \) & \( 0.979 \) & \( 0.954 \) & \( 0.927 \) & \( 0.895 \) & \( 0.870 \) \\
	       & \( 0.009 \) & \( 0.044 \) & \( 0.070 \) & \( 0.085 \) & \( 0.094 \) & \( 0.099 \) \\ \cline{2-7} 
	\( 300 \)  & \( 0.993 \) & \( 0.964 \) & \( 0.935 \) & \( 0.903 \) & \( 0.868 \) & \( 0.836 \) \\
	       & \( 0.008 \) & \( 0.047 \) & \( 0.077 \) & \( 0.108 \) & \( 0.144 \) & \( 0.176 \) \\ \cline{2-7} 
	\( 500 \) & \( 0.996 \) & \( 0.972 \) & \( 0.944 \) & \( 0.917 \) & \( 0.874 \) & \( 0.832 \) \\
	       & \( 0.014 \) & \( 0.053 \) & \( 0.099 \) & \( 0.139 \) & \( 0.166 \) & \( 0.194 \) \\ \cline{2-7} 
	\( 1000 \) & \( 0.990 \) & \( 0.957 \) & \( 0.920 \) & \( 0.882 \) & \( 0.841 \) & \( 0.796 \) \\
	       & \( 0.011 \) & \( 0.050 \) & \( 0.098 \) & \( 0.131 \) & \( 0.171 \) & \( 0.188 \) \\ \cline{2-7} 
	\( 2000 \) & \( 0.991 \) & \( 0.951 \) & \( 0.904 \) & \( 0.850 \) & \( 0.803 \) & \( 0.755 \) \\
	       & \( 0.019 \) & \( 0.048 \) & \( 0.088 \) & \( 0.113 \) & \( 0.124 \) & \( 0.146 \) \\ \cline{2-7} 
	\( 3000 \) & \( 0.994 \) & \( 0.958 \) & \( 0.913 \) & \( 0.863 \) & \( 0.813 \) & \( 0.771 \) \\
	       & \( 0.007 \) & \( 0.033 \) & \( 0.054 \) & \( 0.069 \) & \( 0.079 \) & \( 0.092 \) \\ \hline
	\end{tabular}
\end{table*}

The performance of the proposed procedure is very good except the case
when the sample size is of the same order as the dimension.
However, this regime lies beyond the scope of our results. 
If we have enough data, the method demonstrates very good results even in such challenging situations as recovering a direct sum of several subspaces from non-Gaussian (even not sub-Gaussian) data.

%% file: pp_proofs_march.tex
\newpage
\section{Main proofs} 
\label{Section: Proofs}
This section collects the proofs of the main results.
Some additional technical statements are postponed to the Appendix.

	\subsection{Proof of Theorem \ref{Theorem: BvM projector conjugate}}
The Inverse Wishart prior \( \IWsht_p(\G, p+b-1) \) is conjugate to the multivariate Gaussian distribution, so our pseudo-posterior \( \Ppost{\Su} \) is \\
\( \IWsht_p(\G + n\Se, n+p+b-1) \). We will actively use the following well-known property of the Wishart distribution:
\begin{EQA}
	\Su^{-1} \cond \data 
	& \eqdist &
	\sum_{j=1}^{n+p+b-1} W_j W_j^{\T} \cond \data ,  
\end{EQA}
where \( W_j \cond \data \stackrel{i.i.d.}{\sim} \ND(0, (\G + n\Se)^{-1}) \). 

For shortness in this section we will use the notation \( n_p \eqdef n+p+b-1 \) and we assume that \( b \lesssim p \). 
As we will see, this assumption will help us to simplify the bounds, while the case \( b \gtrsim p \) does not bring any gain.
Moreover, define
\begin{EQA}
	\Su_{n,p} 
	& \eqdef &
	\frac{1}{n_{p}} \G + \frac{n}{n_{p}} \Se
\end{EQA} 
and
\begin{EQA}
	\Eu_{n,p} 
	& \eqdef &
	\frac{1}{n_{p}}\sum_{j=1}^{n_{p}} Z_j Z_j^{\T} - \I_{p} \,,
\label{Enpdef1nj1np}
\end{EQA} 
where \( Z_j \cond \data \stackrel{i.i.d.}{\sim} \ND(0, \I_{p}) \). 
Then \( \Su^{-1} \cond \data \) can be represented as
\begin{EQA}
		\Su^{-1} \cond \data 
		& \eqdist & 
		\Su_{n,p}^{-1/2} \, (\Eu_{n,p} + \I_p) \, \Su_{n,p}^{-1/2}. 
\end{EQA}
We may think that in the ``posterior'' world all randomness comes from \( \Eu_{n,p} \).
Moreover, due to Theorem \ref{Th: Covariance concentration}, (i), there is a random set \( \Upsilon \) such that on this set
\begin{EQA}
	\| \Eu_{n,p} \|_{\infty} 
	&\lesssim &
	\sqrt{\frac{\log(n_p) + p}{n_p}} \leq \sqrt{\frac{\log(n) + p}{n}} \,,
\end{EQA}
and its pseudo-posterior measure
\begin{EQA}
	\Ppost{\Upsilon} 
	& \geq & 
	1 - \frac{1}{n} \, .
\end{EQA}

\paragraph{Step 1}
First, we will need the following lemma.
\begin{lemma} 
\label{Lemma: Posterior contraction}
	The following holds on the random set \( \Upsilon \):
\begin{equation}
	\| \Su - \Se \|_{\infty} 
	\lesssim 
	\| \Se \|_{\infty} \sqrt{\frac{\log(n) + p}{n}} + \frac{\| \G \|_{\infty}}{n}.
\label{SuSe}
\end{equation}
\end{lemma}

\begin{proof}
Since \( \Su^{-1} \cond \data \eqdist \Su_{n,p}^{-1/2} \, (\Eu_{n,p} + \I_p) \, \Su_{n,p}^{-1/2} \), we have
\begin{EQA}
    \Su - \Se 
    &=&
    \Su_{n,p}^{1/2} \, (\Eu_{n,p} + \I_p)^{-1} \, \Su_{n,p}^{1/2} - \Se 
    \\
    &=& 
    \Su_{n,p}^{1/2} \, 
    \left[(\Eu_{n,p} + \I_p)^{-1} - \I_p \right] \, \Su_{n,p}^{1/2} 
    + \Su_{n,p} - \Se.
\end{EQA}
Note that
\begin{EQA}
    \| (\Eu_{n,p} + \I_p)^{-1} - \I_p \|_{\infty} 
    &=& 
    \left\| \sum\limits_{s=1}^{\infty} (-\Eu_{n,p})^s \right\|_{\infty} 
    \\
    & \leq &
    \sum\limits_{s=1}^{\infty} \| \Eu_{n,p} \|_{\infty}^s 
    = 
    \frac{\| \Eu_{n,p} \|_{\infty}}{1-\| \Eu_{n,p} \|_{\infty}} 
    \lesssim \| \Eu_{n,p} \|_{\infty}.
\end{EQA}
Hence,
\begin{EQA}
    \| \Su - \Se \|_{\infty} 
    & \lesssim &
    \| \Su_{n,p} \|_{\infty} \| \Eu_{n,p} \|_{\infty} 
    + \| \Su_{n,p} - \Se \|_{\infty}.
\end{EQA}
Finally, the observations that
\begin{EQA}
    \| \Su_{n,p} \|_{\infty} 
    &\leq &
    \frac{\| \G \|_{\infty}}{n} + \| \Se\|_{\infty}, 
    \\
    \| \Su_{n,p} - \Se \|_{\infty}
    & \leq &
    \frac{\| \G \|_{\infty}}{n} + \frac{n_p - n}{n} \| \Se\|_{\infty},
\end{EQA}
finish the proof.
\end{proof}

The condition on the significant spectral gap for \( \St \) and the bound (\ref{Eq: hat_delta_n}) on the operator norm \( \| \Se - \St \| \) imply a significant spectral gap for the empirical covariance \( \Se \).
The crucial Lemma~\ref{Lemma: Projector decomposition} applied with the central projector \( \Pe \) in place of \( \Pt \) allows to obtain the bound on how close the linear operator
\begin{EQA}
	\LhJ(\Su - \Se)
	& \eqdef &
	\sum_{k \in \IndS} 
	\sum_{l \notin \IndS}
	\frac{\uh_{k} \uh_{k}^{\T}(\Su-\Se)\uh_{l} \uh_{l}^{\T} + \uh_{l} \uh_{l}^{\T}(\Su-\Se)\uh_{k} \uh_{k}^{\T}}
		 {\sigmah_{k} - \sigmah_{l}}
\label{LJSimSihemp}
\end{EQA} 
is to \( \Pu - \Pe \).

\begin{lemma} 
\label{Lemma: Projector concentration Wishart}
	The following holds on the random set \( \Upsilon \):
\begin{EQA}
	\sqrt{n} \| \Pu - \Pe - \LhJ(\Su-\Se)\|_{\Fr} 
	& \lesssim & 
	\erroR_0,
\label{PpostlnPuPr}
\end{EQA}
where
\begin{EQA}
	\erroR_0 
	& \eqdef &
	\sqrt{\frac{\ms_{\subspaceset}}{n}} \,\left( 1 + \frac{\widehat{l}_{\subspaceset}}{\gh_{\subspaceset}}\right)\,
	 \frac{ (\log(n) + p) \| \Se \|_{\infty}^2 + \| \G \|_{\infty}^2/n}{{\gh_{\subspaceset}}^2},
	\label{Delt0defpp}
\end{EQA}
and \( \widehat{l}_{\subspaceset}, \gh_{\subspaceset} \) are empirical versions of \( l^*_{\subspaceset}, \gs_{\subspaceset} \).

\end{lemma}

\begin{proof}
	It follows from (\ref{Formula: Bound}) from Lemma~\ref{Lemma: Projector decomposition} that
\begin{EQA}
	\| \Pu - \Pe - \LhJ(\Su-\Se)\|_{\infty} 
	&\lesssim &
	\left( 1 + \frac{\widehat{l}_{\subspaceset}}{\gh_{\subspaceset}} \right) \, \frac{\| \Su - \Se \|_{\infty}^2}{\gh_{\subspaceset}^2}.
\end{EQA}
It is easy to see that the rank of \( \LhJ(\Su-\Se) \) is at most \( 2 \ms_{\subspaceset} \),
	and thus the rank of \( \Pu - \Pe - \LhJ(\Su-\Se) \) is at most  \( 4 \ms_{\subspaceset} \).
	Hence, taking into account the relation between the Frobenius and the spectral norm of a matrix via rank and 
	(\ref{SuSe}) from Lemma~\ref{Lemma: Posterior contraction}, we obtain the desired statement.
\end{proof}

The representation
\begin{EQA}
		\Su^{-1} \cond \data 
		& \eqdist & 
		\Su_{n,p}^{-1/2} \, (\Eu_{n,p} + \I_p) \, \Su_{n,p}^{-1/2}. 
\end{EQA}
helps to obtain the next result showing that \( \LhJ(\Su-\Se) \) can be approximated by 
\( \StJ = \LhJ \left(-\Se^{1/2} \Eu_{n,p} \Se^{1/2}\right) \).

\begin{lemma}
\label{LLhJSuSe}
It holds
\begin{EQA}
	\LhJ (\Su)
	&=&
	\LhJ \left(-\Se^{1/2} \,\Eu_{n,p}\, \Se^{1/2} \right) + \erroS
	=
	\StJ + \erroS,
\label{LhJSuSEStJerrS}
\end{EQA}
where the remainder \( \erroS \) fulfills on the random set \( \Upsilon \)
\begin{EQA}
	\sqrt{n} \| \erroS \|_{\Fr} 
	\lesssim
	\erroR_{1} 
	& \eqdef &
	 \frac{\ms_{\subspaceset}}{\sqrt{n}} \cdot
	\frac{ (\log(n) + p) \| \Se\|_{\infty} + \| \G \|_{\infty}}{\gh_{\subspaceset}}.
\label{errSR1pp}
\end{EQA}
\end{lemma}

\begin{proof}
Define \( \Ru_{n,p} \) by
\begin{EQA}
	\Ru_{n,p}
	& \eqdef &
	\bigl( \I_{p} + \Eu_{n,p} \bigr)^{-1} - \I_{p} + \Eu_{n,p} . 
\label{RnpdefIpEnpIp}
\end{EQA}
Its spectral norm can be bounded as
\begin{EQA}
	\| \Ru_{n,p} \|_{\infty}
	& \lesssim & 
	\left\| \sum\limits_{s=2}^{\infty} (-\Eu_{n,p})^s \right\|_{\infty} \leq 
     \sum\limits_{s=2}^{\infty} \| \Eu_{n,p} \|_{\infty}^s
     =  
     \frac{\| \Eu_{n,p} \|_{\infty}^2}{1-\| \Eu_{n,p} \|_{\infty}} 
     \lesssim  
     \| \Eu_{n,p} \|_{\infty}^2.
\label{1RnplEnp2lnp}
\end{EQA}
So
\begin{EQA}
	\Su
	& = &
	\Su_{n,p}^{1/2} \, (\Eu_{n,p} + \I_p)^{-1} \, \Su_{n,p}^{1/2}	
	=
	\Su_{n,p}^{1/2} \, (\I_p - \Eu_{n,p} + \Ru_{n,p}) \, \Su_{n,p}^{1/2}.
\end{EQA}
Therefore for \( \Su - \Se \) we have
\begin{EQA}
	\Su - \Se 
	&=&
	\Su_{n,p}^{1/2} \, (\I_p - \Eu_{n,p} + \Ru_{n,p}) \, \Su_{n,p}^{1/2} 
	- \Se 
	\\
	& = &
	-\Su_{n,p}^{1/2} \, \Eu_{n,p} \, \Su_{n,p}^{1/2} 
	+ \Su_{n,p}^{1/2} \, \Ru_{n,p} \, \Su_{n,p}^{1/2} + \Su_{n,p} - \Se.
\end{EQA}
From \( \Su_{n,p}^{1/2} \, \Eu_{n,p} \, \Su_{n,p}^{1/2} \) we pass to \( \Se^{1/2} \, \Eu_{n,p} \, \Se^{1/2} \):
\begin{EQA}
	\Su - \Se
	& = &
	-\Se^{1/2} \, \Eu_{n,p} \, \Se^{1/2} + (\Se^{1/2} \, \Eu_{n,p} \, \Se^{1/2} - \Su_{n,p}^{1/2} \, \Eu_{n,p} \, \Su_{n,p}^{1/2}) \\
	&& \qquad  +   \Su_{n,p}^{1/2} \, \Ru_{n,p} \, \Su_{n,p}^{1/2} + \Su_{n,p} - \Se\\
	& =  &
	-\Se^{1/2} \, \Eu_{n,p} \, \Se^{1/2} + \Ru_1 + \Ru_{2} + \Ru_{3},
\end{EQA}
where we introduce the remainder terms
\begin{EQA}
	\Ru_1 & \eqdef & \Se^{1/2} \, \Eu_{n,p} \, \Se^{1/2} - \Su_{n,p}^{1/2} \, \Eu_{n,p} \, \Su_{n,p}^{1/2} , \\
	\Ru_{2} & \eqdef & \Su_{n,p}^{1/2} \, \Ru_{n,p} \, \Su_{n,p}^{1/2} , \\
	\Ru_{3} & \eqdef & \Su_{n,p} - \Se.
\end{EQA}
They can be bounded as
\begin{EQA}
	\| \Ru_1 \|_{\infty} & \leq & \| \Eu_{n,p} \|_{\infty} \| \Se - \Su_{n,p} \|_{\infty}^{1/2} \left(\| \Su_{n,p}\|_{\infty}^{1/2} + \| \Se \|_{\infty}^{1/2} \right) , \\
	\| \Ru_{2} \|_{\infty} & \leq & \| \Ru_{n,p} \|_{\infty} \| \Su_{n,p} \|_{\infty} \lesssim \| \Eu_{n,p} \|_{\infty}^2 \| \Su_{n,p} \|_{\infty} , \\
	\| \Ru_{3} \|_{\infty} & \lesssim & \frac{\| \G \|_{\infty} + (n_p-n)\,\| \Se \|_{\infty}}{n_p}.
\end{EQA}
Hence, omitting higher order terms, on \( \Upsilon \) we have
\begin{EQA}
	\| \Ru_1 \|_{\infty} & \lesssim &  \| \Se \|_{\infty}^{1/2} \left( \| \G\|_{\infty} + p\| \Se \|_{\infty} \right)^{1/2} \frac{\sqrt{\log(n) + p}}{n}  , \\
	\| \Ru_{2} \|_{\infty} & \lesssim &  \| \Se \|_{\infty} \frac{\log(n) + p}{n}, \\
	\| \Ru_{3} \|_{\infty} & \lesssim & \frac{\| \G \|_{\infty} + p\,\| \Se \|_{\infty}}{n}.
\end{EQA}
Now we summarize
\begin{EQA}
	\LhJ(\Su - \Se) 
	& = & 
	\StJ + \erroS
\label{LJSimSihempex}
\end{EQA} 
with
\begin{EQA}
	\StJ
	& \eqdef &
	-\sum_{k \in \IndS} 
	\sum_{l \not\in \IndS}
	\frac{\sigmah_{k}^{1/2} \, \sigmah_{l}^{1/2} \, (\uh_{k} \uh_{k}^{\T} \Eu_{n,p} \, \uh_{l} \uh_{l}^{\T} + \uh_{l} \uh_{l}^{\T} \Eu_{n,p} \, \uh_{k} \uh_{k}^{\T})}
		 {\sigmah_{k} - \sigmah_{l}},
	\\
	\erroS
	& \eqdef &
	\sum_{k \in \IndS} \sum_{l \not\in \IndS}
	\frac{(\uh_{k} \uh_{k}^{\T} (\Ru_1 + \Ru_{2} + \Ru_{3}) \, \uh_{l} \uh_{l}^{\T} + \uh_{l} \uh_{l}^{\T} (\Ru_1 + \Ru_{2} + \Ru_{3}) \, \uh_{k} \uh_{k}^{\T})}
		 {\sigmah_{k} - \sigmah_{l}} \, .
\label{erroRdef2Rnprst1}
\end{EQA}

Moreover,
\begin{EQA}
	\bigl\| \erroS \bigr\|_{\Fr}
	& \leq &
	2 \left\| \sum_{k \in \IndS} \uh_{k} \uh_{k}^{\T} \sum_{l \not\in \IndS}
	\frac{ (\Ru_1 + \Ru_{2} + \Ru_{3}) \, \uh_{l}\uh_{l}^{\T} }
		 {\sigmah_{k} - \sigmah_{l}} \right\|_2
	\\
	& \leq & 
	2  \sum_{k \in \IndS} \| \uh_{k} \uh_{k}^{\T} \|_{\Fr} \left\| \sum_{l \not\in \IndS}
	\frac{ (\Ru_1 + \Ru_{2} + \Ru_{3}) \, \uh_{l}\uh_{l}^{\T} }
		 {\sigmah_{k} - \sigmah_{l}} \right\|_{\infty}
    \\
	& \leq & 
	2  \sum_{k \in \IndS} \| \Ru_1 + \Ru_{2} + \Ru_{3} \|_{\infty} \left\| \sum_{l \not\in \IndS}
	\frac{  \uh_{l}\uh_{l}^{\T} }
		 {\sigmah_{k} - \sigmah_{l}} \right\|_{\infty}
    \\
    & \leq & 
	 \frac{ 2 \, \ms_{\subspaceset} }{{\gh_{\subspaceset}}} \,\left(\|\Ru_1\|_{\infty} + \| \Ru_2\|_{\infty} + \| \Ru_3 \|_{\infty} \right),
\label{1erro12Rnprstktl}
\end{EQA}
which provides the desired bound.
Similarly, we have
\begin{EQA}
	\bigl\| \StJ \bigr\|_{\Fr}
	& \leq & 
	\frac{ 2 \, \ms_{\subspaceset} \, \| \Se \|_{\infty} }{{\gh_{\subspaceset}}} \, \| \Eu_{n,p} \|_{\infty}
	\lesssim
	\frac{ \ms_{\subspaceset} \, \| \Se \|_{\infty} }{{\gh_{\subspaceset}}} 
	\sqrt{\frac{\log(n) + p}{n}},
\label{StJBound}
\end{EQA}
where the last inequality holds on \( \Upsilon \).

\end{proof}

The results of Lemmas~\ref{Lemma: Projector concentration Wishart} and \ref{LLhJSuSe}
yield on the random set \( \Upsilon \)
\begin{EQA}
	\sqrt{n} \| \Pe - \Pu - \StJ\|_{\Fr} 
	& \lesssim  & 
	\erroR_0 + \erroR_1 .
\end{EQA}
In addition,
\begin{EQA}
	&& \nquad \left| n\| \Pu - \Pe \|_{\Fr}^2 - n\|\StJ\|_{\Fr}^2 \right|
	\\
	& = & 
	n\| \Pu - \Pe - \StJ \|_{\Fr}^2 + 2 \left\langle \sqrt{n} (\Pu - \Pe - \StJ), \sqrt{n} \StJ \right\rangle_{\Fr}  
	\\
	& \leq & n\| \Pu - \Pe - \StJ \|_{\Fr}^2 + 2 \sqrt{n} \| \Pu - \Pe - \StJ \|_{\Fr} \cdot \sqrt{n} \| \StJ \|_{\Fr}.
\end{EQA}
Thus, taking into account the bound for \( \|\StJ\|_{\Fr} \) and neglecting higher order terms, on \( \Upsilon \) we obtain
\begin{equation}
	\left| n\| \Pu - \Pe \|_{\Fr}^2 - n\|\StJ\|_{\Fr}^2 \right| 
	\lesssim 
	\erroR_{2},
\label{BoundDelta2}
\end{equation}
where
\begin{EQA}
	\erroR_{2}
	& \eqdef &
	\left\{ 
		(\log(n) + p) \left( \,\left( 1 + \frac{\widehat{l}_{\subspaceset}}{\gh_{\subspaceset}} \right) \frac{\sqrt{\ms_{\subspaceset}} \| \Se \|_{\infty}}{\gh_{\subspaceset}} + \ms_{\subspaceset}\right) \| \Se \|_{\infty} + \ms_{\subspaceset}\,\| \G \|_{\infty}
	\right\} \times \\
	&& \hspace{7.2cm} \times \frac{\ms_{\subspaceset}\, \|\Se\|_{\infty}}{{\gh_{\subspaceset}}^2} \, \sqrt{\frac{\log(n) + p}{n}}.
\end{EQA}

\paragraph{Step 2}
The norm \( n \| \StJ \|_{\Fr}^2 \) can be decomposed as follows:
\begin{EQA}
	n \| \StJ \|_{\Fr}^2 
	& = &
	2n\sum\limits_{k'=1}^p \sum\limits_{l'=1}^p \sum_{k \in \IndS} \sum_{l \not\in \IndS}
	\frac{\sigmah_k \sigmah_l}{(\sigmah_k - \sigmah_l)^2} \left( \uh_{k'}^{\T} \uh_k\uh_k^{\T} \Eu_{n,p} \uh_l \uh_l^{\T} \uh_{l'} \right)^2 = 
	\\
	& = &
	2n \sum_{k \in \IndS} \sum_{l \not\in \IndS}
	\frac{\sigmah_k \sigmah_l}{(\sigmah_k - \sigmah_l)^2} \left( \uh_k^{\T} \Eu_{n,p} \uh_l \right)^2.
\end{EQA}
Introduce a vector $\xitJ \in \R^{\ms_{\subspaceset} (p - \ms_{\subspaceset})}$ with components
\begin{EQA}
	\xit_{k,l} 
	& = &
	\sqrt{2n}\, \frac{ \sigmah_k^{1/2}\, \sigmah_l^{1/2}}{\sigmah_k - \sigmah_l} \; \uh_k^{\T} \Eu_{n,p} \uh_l,
\end{EQA}
for \( k \in \IndS, l \not\in \IndS \),
ordered in some particular way that will become clear later.
Note that \( n\| \StJ \|_{\Fr}^2 = \| \xitJ \|^2 \) .
Clearly, for each \( k \leq p \) and \( j \leq n_{p} \)
\begin{EQA}
	\eta_{k,j}
	& \eqdef &
	\uh_{k}^{\T}  Z_j \cond \data \stackrel{i.i.d.}{\sim} \ND(0, 1) .
\label{etakjukTZjcXn}
\end{EQA}
Then the components can be rewritten as
\begin{EQA}
	\xit_{k,l} 
	& = &
	\frac{\sqrt{2}\; \sigmah_k^{1/2}\, \sigmah_l^{1/2} }{\sigmah_k - \sigmah_l} \frac{\sqrt{n}}{n_p}
	 \sum\limits_{j=1}^{n_p} \eta_{k,j} \eta_{l,j},
\end{EQA}
for \( k \in \IndS, l \not\in \IndS \).
To understand the covariance structure of \( \xitJ \),
consider one more pair \( (k', l') \) and investigate the covariance:
\begin{EQA}
	\Gammah_{(k,l), (k', l')} 
	& \eqdef &
	\Cov(\xit_{k,l}, \xit_{k', l'} \cond \data)
	\\
	&=&
	\frac{2n}{n_{p}^2} \sum_{j,j'=1}^{n_{p}} 
		\frac{\sigmah_k^{1/2} \, \sigmah_l^{1/2} \, \sigmah_{k'}^{1/2} \, \sigmah_{l'}^{1/2}}
		     {(\sigmah_{k} - \sigmah_{l}) (\sigmah_{k'} - \sigmah_{l'})}
	\E \bigl( \eta_{k,j} \, \eta_{l,j} \, \eta_{k',j'} \, \eta_{l',j'} \cond \data \bigr)
	\\
	&=&
	\frac{2n}{n_p} \delta_{k,k'} \, \delta_{l,l'} \, \frac{\sigmah_{k} \, \sigmah_{l}}
		     {(\sigmah_{k} - \sigmah_{l})^{2}} \, 
\label{2nppkksllstk2}
\end{EQA}
with \( \delta_{k,k'} = \Ind(k=k') \).
Therefore, the covariance matrix of \( \xitJ \) is diagonal:
\begin{EQ}[rcl]
	\GammahJ 
	& \eqdef &
	\frac{2n}{n_p} \cdot \diag\left( 
     	\frac{2 \, \sigmah_{k} \, \sigmah_{l}}{(\sigmah_{k} - \sigmah_{l})^{2}}
    \right)_{k \in \IndS,\, l \notin \IndS} \, .
\label{Formula: tilde_Gamma}
\end{EQ}
This matrix \( \GammahJ \) can be compared with the matrix \( \GammasJ \) defined in (\ref{Formula: Gamma_r}). 
\begin{lemma}
\label{LGamtJsJ}
On the event where \( \| \Se - \St \|_{\infty} \leq \gs_{\subspaceset}/4 \) it holds
\begin{equation}
	\bigl\| \GammahJ - \GammasJ \bigr\|_{1}
	\lesssim 
	\erroR_{3}
\label{GamtJGamsJ1}
\end{equation}
with 
\begin{equation}
	\erroR_{3}
	\eqdef 
	\frac{ p\,\left(\ms_{\subspaceset} \, \| \St \|_{\infty}^2 \land \Tr\left({\St}^2\right) \right) }{{\gs_{\subspaceset}}^3} \left( \| \Se-\St \|_{\infty} + \frac{p}{n} \| \St \|_{\infty} \right).
\nonumber
\end{equation}
\end{lemma}
\begin{proof}
As both matrices \( \GammahJ \) and \( \GammasJ \) are diagonal, it holds
\begin{EQA}
	\bigl\| \GammahJ - \GammasJ \bigr\|_{1}
	&\leq &
	2 \sum_{k \in \IndS}  \sum_{l \notin \IndS}  
		\biggl| 
			\frac{n}{n_p}\frac{\sigmah_{k} \, \sigmah_{l}}{(\sigmah_{k} - \sigmah_{l})^{2}} 
			- \frac{\sigmas_{k} \sigmas_{l}}{(\sigmas_{k} - \sigmas_{l})^{2}} 
		\biggr|.
\label{GtJtGsn1}
\end{EQA}
Let us fix arbitrary \( k \in \IndS, l \notin \IndS\) and upperbound the corresponding term of the sum.
We will extensively use \( |\sigmah_k - \sigmas_k| \leq \| \Se-\St\|_{\infty} \) and \( |\sigmah_l - \sigmas_l| \leq \| \Se-\St\|_{\infty} \) which holds due to the Weyl's inequality.

So, we have
\begin{EQA}
    \biggl| 
			\frac{n}{n_p}\frac{\sigmah_{k} \, \sigmah_{l}}{(\sigmah_{k} - \sigmah_{l})^{2}} 
			- \frac{\sigmas_{k} \sigmas_{l}}{(\sigmas_{k} - \sigmas_{l})^{2}} 
    \biggr|
	& \leq &
	\biggl| 
			\frac{n}{n_p}\frac{\sigmah_{k} \, \sigmah_{l}}{(\sigmah_{k} - \sigmah_{l})^{2}} 
			- \frac{n}{n_p} \frac{\sigmas_{k} \sigmas_{l}}{(\sigmas_{k} - \sigmas_{l})^{2}} 
	\biggr| \\
	&& \quad +
	\biggl| \left( \frac{n}{n_p} - 1 \right) \frac{\sigmas_{k} \sigmas_{l}}{(\sigmas_{k} - \sigmas_{l})^2} \biggr|.
\end{EQA}
Since \( n/n_p \leq 1 \), the first term is controlled by
\begin{EQA}
    && \nquad \biggl| 
		\frac{\sigmah_{k} \, \sigmah_{l}}{(\sigmah_{k} - \sigmah_{l})^{2}} 
		- \frac{\sigmas_{k} \sigmas_{l}}{(\sigmas_{k} - \sigmas_{l})^{2}} 
    \biggr|
	\leq 
    \biggl| 
		\frac{\sigmah_{k} \, \sigmah_{l} \, (\sigmas_{k} - \sigmas_{l})^{2} - \sigmas_{k}\, \sigmas_{l}\, (\sigmah_{k} - \sigmah_{l})^{2} }{(\sigmas_{k} - \sigmas_{l})^{2}\,(\sigmah_{k} - \sigmah_{l})^{2}} 
    \biggr|
    \\
    & = &
    \biggl| 
		\frac{(\sigmas_{k} + \varepsilon_k) \, (\sigmas_{l} + \varepsilon_l) \, (\sigmas_{k} - \sigmas_{l})^{2} - \sigmas_{k}\, \sigmas_{l}\, (\sigmas_{k} - \sigmas_{l} + \varepsilon_k - \varepsilon_l)^{2} }{(\sigmas_{k} - \sigmas_{l})^{2}\,(\sigmah_{k} - \sigmah_{l})^{2}} 
    \biggr|,
\end{EQA}
where we introduced \( \varepsilon_k = \sigmah_{k} - \sigmas_{k} \) and  \( \varepsilon_l = \sigmah_{l} - \sigmas_{l} \). 
Then, crossing out the term \( \sigmas_{k} \, \sigmas_{l}\, (\sigmas_{k} - \sigmas_{l})^{2} \) in the numerator, we obtain
\begin{EQA}
    &&  \biggl| 
		\frac{\sigmah_{k} \, \sigmah_{l}}{(\sigmah_{k} - \sigmah_{l})^{2}} 
		- \frac{\sigmas_{k} \sigmas_{l}}{(\sigmas_{k} - \sigmas_{l})^{2}} 
    \biggr|
	\leq 
    \\
    && \leq
    \biggl| 
		\frac{(\sigmas_{k}\varepsilon_l + \sigmas_{l}\varepsilon_k + \varepsilon_k \varepsilon_l) \, (\sigmas_{k} - \sigmas_{l})^{2} - \sigmas_{k}\, \sigmas_{l}\,( 2(\sigmas_{k} - \sigmas_{l})(\varepsilon_k - \varepsilon_l) + (\varepsilon_k - \varepsilon_l)^{2}) }{(\sigmas_{k} - \sigmas_{l})^{2}\,(\sigmah_{k} - \sigmah_{l})^{2}} 
    \biggr|.
\end{EQA}
On the event where \( \| \Se - \St \|_{\infty} \leq \gs_{\subspaceset}/4 \) we have \( |\varepsilon_k|, |\varepsilon_l| \leq |\sigmas_{k} - \sigmas_{l}|/4\), therefore we can omit the terms \( \varepsilon_k \varepsilon_l\) and \( (\varepsilon_k - \varepsilon_l)^2 \) paying a constant factor for that:
\begin{EQA}
    && \biggl| 
		\frac{\sigmah_{k} \, \sigmah_{l}}{(\sigmah_{k} - \sigmah_{l})^{2}} 
		- \frac{\sigmas_{k} \sigmas_{l}}{(\sigmas_{k} - \sigmas_{l})^{2}} 
    \biggr| \\
	&& \quad \lesssim 
    \biggl| 
		\frac{(\sigmas_{k}\varepsilon_l + \sigmas_{l}\varepsilon_k ) \, (\sigmas_{k} - \sigmas_{l})^{2} - 2\,\sigmas_{k}\, \sigmas_{l}\, (\sigmas_{k} - \sigmas_{l})(\varepsilon_k - \varepsilon_l) }{(\sigmas_{k} - \sigmas_{l})^{2}\,(\sigmah_{k} - \sigmah_{l})^{2}} 
    \biggr| \\
    && \quad = 
    \biggl| 
		\frac{(\sigmas_{k}\varepsilon_l + \sigmas_{l}\varepsilon_k ) \, (\sigmas_{k} - \sigmas_{l}) - 2\,\sigmas_{k}\, \sigmas_{l}\, (\varepsilon_k - \varepsilon_l) }{(\sigmas_{k} - \sigmas_{l}) \,(\sigmah_{k} - \sigmah_{l})^{2}} 
    \biggr| \\
    && \quad =
    \biggl| 
		\frac{\varepsilon_k ({\sigmas_{k}}^2 + \sigmas_{k} \sigmas_{l}) - \varepsilon_l ({\sigmas_{l}}^2 + \sigmas_{k} \sigmas_{l})}{(\sigmas_{k} - \sigmas_{l}) \,(\sigmah_{k} - \sigmah_{l})^{2}} 
    \biggr|\, .
\end{EQA}
Since \( \sigmas_{k} \sigmas_{l} \leq ({\sigmas_{k}}^2 + {\sigmas_{l}}^2)/2\), we get
\begin{EQA}
    \biggl| 
		\frac{\sigmah_{k} \, \sigmah_{l}}{(\sigmah_{k} - \sigmah_{l})^{2}} 
		- \frac{\sigmas_{k} \sigmas_{l}}{(\sigmas_{k} - \sigmas_{l})^{2}} 
    \biggr| 
	& \lesssim &
		\frac{{\sigmas_{k}}^2 + {\sigmas_{l}}^2}{(\sigmas_{k} - \sigmas_{l}) \,(\sigmah_{k} - \sigmah_{l})^{2}}\, \| \Se-\St \|_{\infty} .
\end{EQA}
As to the denominator, again considering the event where \( \| \Se - \St \|_{\infty} \leq \gs_{\subspaceset}/4 \), we have \( |\varepsilon_k|, |\varepsilon_l| \leq |\sigmas_{k} - \sigmas_{l}|/4\) that results in
\begin{EQA}
      |\sigmah_{k} - \sigmah_{l}| 
      & = &
      |\sigmas_{k} - \sigmas_{l} + \varepsilon_k - \varepsilon_l|
      \geq
       |\sigmas_{k} - \sigmas_{l}|/2.
\end{EQA}
Hence,
\begin{EQA}
    \biggl| 
		\frac{\sigmah_{k} \, \sigmah_{l}}{(\sigmah_{k} - \sigmah_{l})^{2}} 
		- \frac{\sigmas_{k} \sigmas_{l}}{(\sigmas_{k} - \sigmas_{l})^{2}} 
    \biggr| 
	& \lesssim &
		\frac{{\sigmas_{k}}^2 + {\sigmas_{l}}^2}{(\sigmas_{k} - \sigmas_{l})^3} \, \| \Se-\St \|_{\infty} \,.
\end{EQA}

As to the term \( \biggl| \left( \frac{n}{n_p} - 1 \right) \frac{\sigmas_{k} \sigmas_{l}}{(\sigmas_{k} - \sigmas_{l})^2} \biggr| \), it is simply bounded as
\begin{EQA}
      \biggl| \left( \frac{n}{n_p} - 1 \right) \frac{\sigmas_{k} \sigmas_{l}}{(\sigmas_{k} - \sigmas_{l})^2} \biggr|
      & \lesssim &
      \frac{p}{n} \cdot \frac{{\sigmas_{k}}^2 + {\sigmas_{l}}^2}{(\sigmas_{k} - \sigmas_{l})^2}
      \leq
      \frac{p}{n} \cdot \frac{{\sigmas_{k}}^2 + {\sigmas_{l}}^2}{(\sigmas_{k} - \sigmas_{l})^3} \, \| \St \|_{\infty} \,.
\end{EQA}
The last inequality uses the fact that \( |\sigmas_{k} - \sigmas_{l} | \leq \| \St \|_{\infty} \) which is rather useless, however it allows to write the final bound in a convenient form and doesn't worsen the result.

Putting this all together, we get
\begin{EQA}
	\bigl\| \GammahJ - \GammasJ \bigr\|_{1} 
	& \lesssim &
	\sum_{k \in \IndS} \sum_{l \notin \IndS} 
	\frac{{\sigmas_k}^2 + {\sigmas_l}^2 }{(\sigmas_k-\sigmas_l)^3} \left( \| \Se-\St \|_{\infty} + \frac{p}{n} \, \| \St \|_{\infty} \right),
\end{EQA}
which provides the desired result once we notice that \( \sum\limits_{k \in \IndS} \sum\limits_{l \notin \IndS}
	({\sigmas_k}^2 + {\sigmas_l}^2) \) can be bounded by both \( 2\,p\, \Tr ({\St}^2)\) and \( 2\,p \, \ms_{\subspaceset} \, \| \St \|_{\infty} \).
\end{proof}

Unfortunately, the entries \( \xit_{k,l} \) of \( \xitJ \) are not Gaussian because of the product \( \eta_{k,j} \, \eta_{l,j} \).
This does not allow to apply the Gaussian comparison Lemma~\ref{Lemma: Gaussian comparison}.
To get rid of this issue, we condition on \( \Pe Z \). Namely, 
in the ``posterior'' world random vectors \( \Pe Z_j \) and \( (\I_p -\Pe) Z_j \) are Gaussian and uncorrelated, therefore, independent, so we can condition on \( Z_{\subspaceset} \eqdef (\Pe Z_{1}, \ldots, \Pe Z_{n_{p}}) \) to get that \( \StJ \) is conditionally on \( \data, Z_{\subspaceset} \) Gaussian random vector with the covariance matrix
\begin{equation}
\begin{aligned}
	\GammatJ
	& \eqdef &
	\Cov\bigl( \xitJ \cond \data, Z_{\subspaceset} \bigr) .
\nonumber
\end{aligned}
\end{equation}
It holds similarly to the above
\begin{equation}
\begin{aligned}
	\Gammat_{(k,l), (k', l')} 
	& \eqdef 
	\Cov(\xit_{k,l}, \xit_{k', l'} \cond \data, Z_{\subspaceset})
	\\
	&=
	\frac{2n}{n_{p}^2} \sum_{j,j'=1}^{n_{p}} 
		\frac{\sigmah_{k}^{1/2} \, \sigmah_{l}^{1/2} \, \sigmah_{k'}^{1/2} \, \sigmah_{l'}^{1/2}}
		     {(\sigmah_{k} - \sigmah_{l}) (\sigmah_{k'} - \sigmah_{l'})}
	\E \bigl( 
	 \eta_{k,j} \, \eta_{l,j} \, \eta_{k',j'} \, \eta_{l',j'} \cond \data, Z_{\subspaceset} 
	\bigr)
	\\
	&=
	\frac{2n}{n_p} \; \tdelta_{k,k'} \delta_{l,l'} \frac{\sigmah_{k}^{1/2} \, \sigmah_{k'}^{1/2} \, \sigmah_{l}}
		     {(\sigmah_{k} - \sigmah_{l}) (\sigmah_{k'} - \sigmah_{l})} 
\nonumber
\end{aligned}
\end{equation}
with
\begin{EQA}
	\tdelta_{k,k'}
	& \eqdef &
	\frac{1}{n_{p}} \sum_{j=1}^{n_{p}} \eta_{k,j} \eta_{k',j} \, .
\label{tdkks1npj1}
\end{EQA}

\begin{lemma}
\label{GammatGammahpp}
It holds on a random set of pseudo-posterior measure \( 1 - \frac{1}{n} \)
\begin{EQA}
	\max_{k,k' \in \IndS} \bigl| \tdelta_{k,k'} - \delta_{k,k'} \bigr|
	& \lesssim &
	\sqrt{\frac{\log(n_{p} + \ms_{\subspaceset})}{n_{p}}} \,,
\label{maxkkaIStdkka2lnp}
\end{EQA}
and on this set
\begin{equation}
	\bigl\| \GammatJ - \GammahJ \bigr\|_{1}
	\lesssim 
	\erroR_{4}
	\eqdef
	\frac{ (\ms_{\subspaceset})^{3/2} \| \Se \|_{\infty} \, \Tr(\Se)}{\gh_{\subspaceset}^{2}}
	\sqrt{\frac{\log(n_{p}+\ms_{\subspaceset})}{n_{p}}} \, .
\label{GtGh1e442np}
\end{equation}
\end{lemma}
\begin{proof}
	The first result of the lemma follows easily from usual concentration inequalities for sub-exponential random variables and union bound for at most \( |\IndS|^2 = (\ms_{\subspaceset})^2\) pairs of \( k, k' \).

	To obtain the second inequality we represent \( \GammatJ\) and \( \GammahJ \) as
\begin{EQA}
		\GammatJ & = & \diag\left( \GammatJ^{(l)}\right)_{l \notin \IndS},\\
		\GammahJ & = & \diag\left( \GammahJ^{(l)}\right)_{l \notin \IndS}.
\end{EQA}
	Due to this block structure we have
\begin{EQA}
		\| \GammatJ - \GammahJ \|_1 
		& = & 
		\sum\limits_{l \notin \IndS} \| \GammatJ^{(l)} - \GammahJ^{(l)} \|_1.
\end{EQA}
	Let us fix \( l \notin \IndS \) and focus on the corresponding block with size \( \ms_{\subspaceset} \times \ms_{\subspaceset} \).
	It's easy to observe that for each \( k, k' \in \IndS \)
\begin{EQA}
		\Gammat_{(k,l), (k', l)} - \Gammah_{(k,l), (k', l)}
		& = &
		\frac{2n}{n_p}\; 
		 \frac{\sigmah_k^{1/2} \, \sigmah_{k'}^{1/2} \, \sigmah_l}{(\sigmah_k-\sigmah_l)(\sigmah_{k'} - \sigmah_l)}  
		\cdot ( \tdelta_{k,k'} - \delta_{k,k'} )
\end{EQA}
	and, therefore,
\begin{EQA}
		\max\limits_{k,k' \in \IndS} \left|\Gammat_{(k,l), (k', l)} - \Gammah_{(k,l), (k', l)} \right|
		& \leq &
		\frac{2 \| \Se \|_{\infty} \sigmah_l}{\gh_{\subspaceset}^2}\;
		\max\limits_{k,k' \in \IndS} \bigl| \tdelta_{k,k'} - \delta_{k,k'} \bigr|.
\end{EQA}
	Finally, since
\begin{EQA}
		\| \GammatJ^{(l)} - \GammahJ^{(l)} \|_1 
		& \leq & 
		\sqrt{\ms_{\subspaceset}} \| \GammatJ^{(l)} - \GammahJ^{(l)} \|_{2} \\
		& \leq & 
		(\ms_{\subspaceset})^{3/2} \max\limits_{k,k' \in \IndS} \left|\Gammat_{(k,l), (k', l')} - \Gammah_{(k,l), (k', l')} \right|,
\end{EQA}{}
	the obtained inequalities provide the result of the lemma.
\end{proof}

Putting together (\ref{GamtJGamsJ1}) and (\ref{GtGh1e442np})  yields the bound
\begin{EQA}
	\bigl\| \GammatJ - \GammasJ \bigr\|_{1}
	& \lesssim &
	\erroR_{3} + \erroR_{4} .
\label{GJtGJs1e3e4}
\end{EQA}
The Gaussian comparison Lemma~\ref{Lemma: Gaussian comparison} can be used to compare the conditional distribution
of \( \| \xitJ \| \) given \( \data, \Pe Z \) and the unconditional distribution of \( \| \xiJ \| \):
on a random set of pseudo-posterior measure \( 1 - \frac{1}{n} \)
\begin{EQA}
	&& \sup_{x\in\R} \left| 
		\prior \left( \| \xitJ \|^2 \leq x \Cond \data, Z_{\subspaceset} \right) 
			- \Pro\left( \|\xiJ\|^2 \leq x\right) 
	\right| \\
	&& \quad \lesssim 
	\frac{\erroR_{3} + \erroR_{4}}{\| \Gammas_{\subspaceset} \|_2^{1/2} \left( \| \Gammas_{\subspaceset} \|_2^2 - \| \Gammas_{\subspaceset} \|_{\infty}^2 \right)^{1/4}} \, .
\end{EQA} 
Of course, integrating  w.r.t. \( \Pe Z \) ensures similar result when conditioning on the data \( \data \) only:
\begin{equation}
\begin{aligned}
	& \sup_{x\in\R} \left| 
		\prior \left( \| \xitJ \|^2 \leq x \Cond \data \right) 
			- \Pro\left( \|\xiJ\|^2 \leq x\right) 
	\right| \\
	& \quad \lesssim 
	\frac{\erroR_{3} 
	+ \erroR_{4}}{\| \Gammas_{\subspaceset} \|_2^{1/2} \left( \| \Gammas_{\subspaceset} \|_2^2 - \| \Gammas_{\subspaceset} \|_{\infty}^2 \right)^{1/4}} 
	+ \frac{1}{n}
\label{Bound: delta34}
\end{aligned}
\end{equation} 
with probability one.

\paragraph{Step 3 }
		So far our bounds \( \widehat{\Delta}_2, \widehat{\Delta}_3, \widehat{\Delta}_4 \) were obtained in the ``posterior'' world, and they are random in the \(\data-\)world since they depend on \( \Se \).
		We want to bound them by deterministic counterparts \( \Delta_1, \Delta_2, \Delta_3\) with probability \( 1 - 1/n \).
		To do so, we basically need to upperbound \( \| \Se \|_{\infty}, \; \Tr(\Se), \; \widehat{l}_{\subspaceset} \) and to lowerbound \( \gh_{\subspaceset} \).
		This can be done as follows:
		\begin{EQA}
			&& \| \Se \|_{\infty} \leq \| \St \|_{\infty} + \| \Se - \St \|_{\infty} \leq \| \St \|_{\infty} (1 + \hdelta_n),\\
			&& \Tr(\Se) \leq \Tr(\St) + p\| \Se-\St\|_{\infty} \leq \Tr(\St) \left( 1 + \frac{p\,\hdelta_n\,\| \St \|_{\infty}}{\Tr(\St)}\right).
		\end{EQA}
		with probability \( 1 - 1/n \). Due to the definition of the spectral gap, there is an index \( j \) such that \( \gh_{\subspaceset} = \widehat{\sigma}_j - \widehat{\sigma}_{j+1}\). Then due to the Weyl's inequality we have
			\begin{EQA}
			\gh_{\subspaceset} & = & \widehat{\sigma}_j - \widehat{\sigma}_{j+1} \geq 
			\sigma^*_j - \sigma^*_{j+1} - |\widehat{\sigma}_j - \sigma^*_j| - |\widehat{\sigma}_{j+1} - \sigma^*_{j+1}| \\
			& \geq &
			\gs_{\subspaceset} - 2 \| \Se - \St \|_{\infty} \geq
			\gs_{\subspaceset} - 2\hdelta_n \| \St \|_{\infty}
		\end{EQA}
		with probability \( 1 - 1/n \). Similarly we can upperbound \( \widehat{l}_{\subspaceset}\).
		Further we can plug the obtained inequalities directly in \( \widehat{\Delta}_2, \widehat{\Delta}_3, \widehat{\Delta}_4 \) in order to get \( \Delta_1, \Delta_2, \Delta_3\) without any assumptions on \( \hdelta_n \).
		Another option is to use our assumption 
		\begin{EQA}
			\hdelta_n & \leq & \frac{\gs_{\subspaceset}}{4\| \St \|_{\infty}} \;\land\; \frac{r(\St)}{p},
		\end{EQA}
		that ensures 
		\( \| \Se \|_{\infty} \lesssim \| \St \|_{\infty},\;
		 \Tr(\Se) \lesssim \Tr(\St),\;
		  \gh_{\subspaceset} \gtrsim \gs_{\subspaceset},\; 
		  \widehat{l}_{\subspaceset} \lesssim l^*_{\subspaceset} \).
		This allows to obtain  more transparent bounds on \( \Delta_1, \Delta_2, \Delta_3\). 
		Also note that this assumption guarantees that the event from Lemma~\ref{LGamtJsJ} is of probability at least \( 1 - 1/n \).
		We conclude that
\begin{EQA}
	\erroR_{2} 
	& \lesssim & 
	\erroD_{2} 
	\eqdef
	\left\{ 
		(\log(n) + p) \left( \,\left( 1 + \frac{l^*_{\subspaceset}}{\gs_{\subspaceset}} \right) \frac{\sqrt{\ms_{\subspaceset}} \| \St \|_{\infty}}{\gs_{\subspaceset}} + 1 \right) \| \St \|_{\infty} + \| \G \|_{\infty}
	\right\} \times \\
	&& \hspace{7.4cm} \times \frac{\ms_{\subspaceset} \, \| \St\|_{\infty}}{{\gs_{\subspaceset}}^2} \sqrt{\frac{\log(n) + p}{n}}
	,\\
	\erroR_{3}
	& \lesssim & 
	\erroD_{3}
	\eqdef
	\frac{ \| \St \|_{\infty} \,\left(\ms_{\subspaceset} \, \| \St \|_{\infty}^2 \land \Tr\left({\St}^2\right) \right) }{{\gs_{\subspaceset}}^3} \, p\left( \hdelta_n + \frac{p}{n} \right),\\
	\erroR_4 
	& \lesssim &
	\erroD_4
	\eqdef 
	\frac{ (\ms_{\subspaceset})^{3/2} \| \St \|_{\infty} \, \Tr(\St)}{{\gs_{\subspaceset}}^{2}}
	\sqrt{\frac{\log(n)}{n}}
\end{EQA}
with probability \( 1 - {1}/{n} \) in the $\data-$world.

Now we combine the obtained bounds.
For \( \erroD_{2} \) defined above and arbitrary \( x \in \R \) it holds
\begin{EQA}
	&& \nquad \Ppost{ n \| \Pu - \Pe \|_{\Fr}^2 \leq x} 
    \\
    && \leq 
    \Ppost{ n \| \StJ\|_{\Fr}^2 \leq x + \erroD_{2}}  
    \\
    && \qquad
    + \, \Ppost{ n \| \Pu - \Pe \|_{\Fr}^2 - n \| \StJ\|_{\Fr}^2 \leq - \erroD_{2}}.
\label{PpPuPrFrgx01}
\end{EQA}
Since \( n \| \StJ\|_{\Fr}^2 \cond \data \eqdist \| \xitJ \|^2 \cond \data \), \( \erroR_{2} \lesssim \erroD_{2}\) with probability \( 1 - \frac{1}{n} \), 
and taking (\ref{BoundDelta2}) into account, we deduce
\begin{EQA}
	\Ppost{ n \| \Pu - \Pe \|_{\Fr}^2 \leq x} 
    & \leq &
    \Ppost{ \| \xiJ\|^2 \leq x + \erroD_{2}} + \Ppost{ \Upsilon^c}
\label{PpPuPrFrgx01}
\end{EQA}
with probability \( 1 - \frac{1}{n} \). Subtracting \( \Pro\left( \| \xiJ\|^2 \leq x  \right) \) and taking supremum of both sides, we get
\begin{EQA}
	&& \nquad 
	\sup_{x\in\R} \left\{ 
		\Ppost{ n \| \Pu - \Pe \|_{\Fr}^2 \leq x} 
		- \Pro\left( \| \xiJ\|^2 \leq x  \right) 
	\right\} 
    \\
    & \leq &
    \sup_{x\in\R} \left\{ \Ppost{ \| \xitJ\|^2 \leq x + \erroD_{2}} - \Pro\left( \| \xiJ\|^2 \leq x+\erroD_{2}  \right) \right\} 
    \\
    && \, + 
    \sup_{x\in\R} 
    \left\{ \Pro\left( \| \xiJ\|^2 \leq x+\erroD_{2}  \right) - \Pro\left( \| \xiJ\|^2 \leq x  \right) \right\} +
    \Ppost{ \Upsilon^c}.
\label{PpPuPrFrgx01}
\end{EQA}
The first term in the right-hand side is bounded by \( \frac{\erroD_{3} + \erroD_{4}}{\| \Gammas_{\subspaceset} \|_2^{1/2} \left( \| \Gammas_{\subspaceset} \|_2^2 - \| \Gammas_{\subspaceset} \|_{\infty}^2 \right)^{1/4}} + \frac{1}{n} \) with probability \( 1 - \frac{1}{n} \) due to (\ref{Bound: delta34}).
The second term does not exceed \( \frac{\erroD_{2}}{\| \Gammas_{\subspaceset} \|_2^{1/2} \left( \| \Gammas_{\subspaceset} \|_2^2 - \| \Gammas_{\subspaceset} \|_{\infty}^2 \right)^{1/4}} \) according to the Gaussian anti-concentration  Lemma~\ref{Lemma: Anticoncentration}. The last term is at most $\frac{1}{n}$ by definition of \( \Upsilon \).
Therefore,
\begin{EQA}
	&& \sup_{x\in\R} \left\{ \Ppost{ n \| \Pu - \Pe \|_{\Fr}^2 \leq x} - \Pro\left( \| \xiJ\|^2 \leq x  \right) \right\} \\
    && \quad \lesssim 
    \frac{\erroD_{2} + \erroD_{3} + \erroD_{4}}{\| \Gammas_{\subspaceset} \|_2^{1/2} \left( \| \Gammas_{\subspaceset} \|_2^2 - \| \Gammas_{\subspaceset} \|_{\infty}^2 \right)^{1/4}} +
    \frac{1}{n}
\end{EQA}
with probability \( 1 - \frac{1}{n}\).
Similarly, one derives
\begin{EQA}
	&& \sup_{x\in\R} \left\{ \Pro\left( \| \xiJ\|^2 \leq x  \right) - \Ppost{ n \| \Pu - \Pe \|_{\Fr}^2 \leq x} \right\} \\
    && \quad \lesssim 
    \frac{\erroD_{2} + \erroD_{3} + \erroD_{4}}{\| \Gammas_{\subspaceset} \|_2^{1/2} \left( \| \Gammas_{\subspaceset} \|_2^2 - \| \Gammas_{\subspaceset} \|_{\infty}^2 \right)^{1/4}} +
    \frac{1}{n}
\end{EQA}
with probability \( 1 - \frac{1}{n}\). The previous two inequalities yield the desired result.

	\subsection{Proof of Corollary \ref{Corollary: Confidence set Wishart}}
Let \( \xiJ \sim \ND(0, \GammasJ) \).
Due to Theorem \ref{Theorem: clt} we have
\begin{EQA}
					&\sup_{x \in \R} 
					\left|
					 \Pro \left( n \| \Pe - \Pt \|_{\Fr}^2  >  x \right) -
					 \Pro(\| \xiJ \|^2 > x) 
					\right| \lesssim \erro.
\end{EQA}
			Fix arbitrary significance level \( \alpha \in (0;\;1) \) (or confidence level \( 1-\alpha \)). 
			Recall that by \( \gamma_{\alpha} \) we denote \( \alpha \)-quantile of \( n \| \Pe - \Pt \|_{\Fr}^2 \).
			Let us fix an event \( \Theta \) such that
\begin{EQA}
					\sup_{x \in\R}\left|
					 \Ppost{n \| \Pu - \Pe \|_{\Fr}^2  >  x} -
					 \Pro(\| \xiJ \|^2 > x) 
					\right| 
					&\lesssim &
					{\Diamond}.
\end{EQA}
			According to Theorem \ref{Theorem: BvM projector conjugate} its probability is at least \( 1 - {1}/{n} \).
			Hence, by the triangle inequality it holds on \( \Theta \)
\begin{EQA}
					\sup_{x\in\R}\left|
					\Ppost{n \| \Pu - \Pe \|_{\Fr}^2  >  x} - 
					\Pro \left( n \| \Pe - \Pt \|_{\Fr}^2  >  x \right)
					\right| 
					&\leq &
					\Diamond^{\prime} \asymp \erro + \Diamond.
\end{EQA}
			Therefore, taking \( x = \gamma_{\alpha-\Diamond^{\prime}} \) and \( x = \gamma_{\alpha+\Diamond^{\prime}} \), we get on \( \Theta \)
\begin{EQA}
					& \left|
					\Ppost{ n \| \Pu - \Pe \|_{\Fr}^2 > \gamma_{\alpha-\Diamond^{\prime}} } - (\alpha - \Diamond^{\prime})
					\right| 
					\leq \Diamond^{\prime},\\
					& \left|
					\Ppost{  n \| \Pu - \Pe \|_{\Fr}^2 > \gamma_{\alpha+\Diamond^{\prime}} } - (\alpha + \Diamond^{\prime})
					\right| 
					\leq \Diamond^{\prime}.
\end{EQA}
			Thus,
\begin{EQA}
					& \Ppost{ n \| \Pu - \Pe \|_{\Fr}^2 > \gamma_{\alpha-\Diamond^{\prime}} }
					\leq (\alpha - \Diamond^{\prime}) + \Diamond^{\prime} = \alpha,\\
					& \Ppost{ n \| \Pu - \Pe \|_{\Fr}^2 > \gamma_{\alpha+\Diamond^{\prime}} }
					\geq (\alpha + \Diamond^{\prime}) - \Diamond^{\prime} = \alpha.
\end{EQA}
By definition of \( \gamma^{\circ}_{\alpha} \) the previous two inequalities yield
\begin{EQA}
	\gamma_{\alpha + \Diamond^{\prime}} \leq \gamma_{\alpha}^{\circ} 
	&\leq &
	\gamma_{\alpha - \Diamond^{\prime}} \;\;\; \text{on }\;\Theta.
\end{EQA}
Hence,
\begin{EQA}
	\Pro \left( \gamma^{\circ}_{\alpha} < \gamma_{\alpha+\Diamond^{\prime}} \right) 
	&\leq &
	\Pro\left({\Theta}^c\right) \leq \frac{1}{n},
	\\
	\Pro \left( \gamma^{\circ}_{\alpha} > \gamma_{\alpha-\Diamond^{\prime}} \right)
	& \leq &
	\Pro\left({\Theta}^c\right) 
	\leq \frac{1}{n}.
\end{EQA}
Now we can write the following chain of inequalities:
\begin{EQA}
	&& \nquad 
	\Pro \left( n \| \Pe - \Pt \|_{\Fr}^2 > \gamma_{\alpha}^{\circ} \right) 
	\\
	& \leq &
	\Pro \left( 
		\left\{ n \| \Pe - \Pt \|_{\Fr}^2 > \gamma_{\alpha+\Diamond^{\prime}} 
		\right\} \cup 
		\left\{ \gamma_{\alpha}^{\circ} < \gamma_{\alpha+\Diamond^{\prime} } 
		\right\} 
	\right) 
	\\
	& \leq &
	\Pro \left( n \| \Pe - \Pt \|_{\Fr}^2 > \gamma_{\alpha+\Diamond^{\prime}} \right) 
	+ \Pro \left( \gamma_{\alpha}^{\circ} < \gamma_{\alpha+\Diamond^{\prime}} \right) 
	\leq 
	\alpha + \Diamond^{\prime} + \frac{1}{n}
\end{EQA}
and
\begin{EQA}
	&& \nquad 
	\Pro \left( n\| \Pe - \Pt \|_{\Fr}^2 > \gamma_{\alpha}^{\circ} \right) 
	=
	1 - \Pro \left( n\| \Pe - \Pt \|_{\Fr}^2 \leq \gamma_{\alpha}^{\circ} \right) 
	\\
	& \geq &
	1 - \Pro \left( 
		\left\{ n\| \Pe - \Pt \|_{\Fr}^2 \leq \gamma_{\alpha-\Diamond^{\prime}} 
		\right\} \cup \left\{ \gamma_{\alpha}^{\circ} > \gamma_{\alpha-\Diamond^{\prime}} 
		\right\} 
	\right) 
	\\
	& \geq &
	1 - \Pro \left( n\| \Pe - \Pt \|_{\Fr}^2 \leq \gamma_{\alpha-\Diamond^{\prime}} \right) 
	- \Pro \left( \gamma_{\alpha}^{\circ} > \gamma_{\alpha-\Diamond^{\prime}} \right) 
	\\
	& = &
	\Pro \left( n\| \Pe - \Pt \|_{\Fr}^2 > \gamma_{\alpha-\Diamond^{\prime}} \right) - \Pro \left( \gamma_{\alpha}^{\circ} > \gamma_{\alpha-\Diamond^{\prime}} \right) 
	\geq
	\alpha - \Diamond^{\prime} - \frac{1}{n}.
\end{EQA}
Finally, these inequalities imply the following bound
\begin{EQA}
	\left| 
		\alpha 
		- \Pro \left( n\| \Pe - \Pt \|_{\Fr}^2 > \gamma_{\alpha}^{\circ} \right) 
	\right| 
	&\leq &
	\Diamond^{\prime} + \frac{1}{n},
\end{EQA}
which concludes the proof.
  
  \section{Acknowledgements}
    We thank an anonymous Referee for very valuable comments and suggestions
which led to a significant improvement of the paper.
    
    This work has been funded by the Russian Academic Excellence Project `5-100'.
    Results of Section 2 has been obtained under support of the RSF grant No. 18-11-00132. Financial support by the German Research Foundation (DFG) through the Collaborative Research Center 1294 is gratefully acknowledged.

	\appendix
	\section{Auxiliary results} \label{Appendix: A}
Here we formulate some well-known results that were used throughout the paper. 

The following theorem gathers several crucial results on concentration of sample covariance.
\begin{theorem} \label{Th: Covariance concentration}
	Let \( X_1, \ldots, X_n \) be i.i.d. zero-mean random vectors in \( \R^p \). 
	Denote the true covariance matrix as \( \St \eqdef \E \left(X_i X_i^{\T}\right) \) and the sample covariance as \( \Se \eqdef \frac{1}{n}\sum\limits_{i=1}^n X_i X_i^{\T} \).
	Suppose the data are obtained from:\\
	(i) Gaussian distribution \( \ND(0, \St) \). In this case, define \( \hdelta_n \) as
\begin{EQA}
	\hdelta_n 
	&\asymp & 
	\sqrt{\frac{r(\St) + \log(n)}{n}} ;
\end{EQA}
	(ii) Sub-Gaussian distribution. In this case, define \( \hdelta_n \) as
\begin{EQA}
	\hdelta_n
	& \asymp &
	\sqrt{\frac{p+\log(n)}{n}};
\end{EQA}
	(iii) a distribution supported in some centered Euclidean ball of radius \( R \). In this case, define \( \hdelta_n \) as
\begin{EQA}
	\hdelta_n 
	&\asymp & 
	\frac{R}{\sqrt{\| \St \|}} \sqrt{\frac{\log(n)}{n}};
\end{EQA}
	(iv) log-concave probability measure. In this case, define \( \hdelta_n \) as
\begin{EQA}
	\hdelta_n
	& \asymp & 
	\sqrt{\frac{\log^6(n)}{np}} .
\end{EQA}
	
	Then in all the cases above the following concentration result for \( \Se \) holds with the corresponding \( \hdelta_n \):
\begin{EQA}
	\|\Se - \St\|_{\infty} 
	&\leq &
	\hdelta_n \| \St\|_{\infty}
\end{EQA}
	with probability at least \( 1-\frac{1}{n} \).
\end{theorem}
\begin{proof}
	(i) See \cite{Koltchinskii_CIAMBFSCO}, Corollary 2.
	(ii) This is a well-known simple result presented in a range of papers and lecture notes. See, e.g. \cite{Rigollet}, Theorem 4.6.
	(iii) See \cite{Vershynin_ITTNAAORM}, Corollary 5.52. Usually the radius \( R \) is taken such that \( \frac{R}{\sqrt{\| \St \|}} \asymp \frac{\sqrt{\Tr(\St)}}{\sqrt{\| \St \|}} = \sqrt{r(\St)} \).
	(iv) See \cite{Adamczak}, Theorem 4.1.
\end{proof}

The following lemma is a crucial tool when working with spectral projectors.
\begin{lemma} \label{Lemma: Projector decomposition}
	The following bound holds for all \( \subspaceset = \{ r^-, r^- + 1, \ldots, r^+ \} \) with \( 1 \leq r^- \leq r^+ \leq q\):
\begin{EQA}
	\| \Pp - \Pt \|_{\infty}
	& \leq &
	4 \left( 1 + \frac{2}{\pi} \frac{l^*_{\subspaceset}}{\gs_{\subspaceset}} \right) \frac{\| \Sp - \St \|_{\infty}}{\gs_{\subspaceset}}.
\label{Eq: Davis-Kahan}
\end{EQA}
	Moreover, the following representation holds:
\begin{equation}
	\Pp - \Pt 
	=
	L_{\subspaceset}(\Sp-\St) + R_{\subspaceset}(\Sp-\St),
\label{Formula: Projector decomposition}
\end{equation}
	where
\begin{EQA}
	L_{\subspaceset}(\Sp-\St)
	& \eqdef &
	\sum_{r \in {\subspaceset}} \sum_{s \notin {\subspaceset}} 
			\frac{\Ptr (\Sp-\St) \Pts + \Pts (\Sp-\St) \Ptr}{\mus_r - \mus_s}
\end{EQA}
	and
\begin{equation}
	\| R_{\subspaceset}(\Sp-\St) \|_{\infty} 
	\leq 
	15\, \left( 1 + \frac{2}{\pi} \frac{l^*_{\subspaceset}}{\gs_{\subspaceset}} \right) \,\left( \frac{\| \Sp-\St \|_{\infty}}{\gs_{\subspaceset}} \right)^2.
\label{Formula: Bound}
\end{equation}
\end{lemma}
\begin{proof}
Apply Lemma 2 from \cite{Koltchinskii_AACBFBFOSPOSC}.
\end{proof}

\noindent 
This lemma shows that \( \Pp - \Pt \) can be approximated by the linear operator \( L_{\subspaceset}(\Sp-\St) \). 

The next lemma from \cite{Goetze} provides upper bound for \( \Delta- \)band of the squared norm of a Gaussian element.
\begin{lemma} [Gaussian anti-concentration]\label{Lemma: Anticoncentration}	
	Let \( \xi \) be a Gaussian element in Hilbert space \( \mathbb{H} \) with zero mean and covariance operator \( \Su_{\xi} \).
	Then for arbitrary \( \Delta > 0 \) one has
\begin{EQA}
	\Pro(x < \| \xi \|^2 < x+\Delta)
	& \leq &
	\frac{\Delta}{\| \Su_{\xi}\|_2^{1/2} \left(\| \Su_{\xi}\|_2^2 - \|\Su_{\xi}\|_{\infty}^2\right)^{1/4}}.
\label{F: Anticoncentration}
\end{EQA}
\end{lemma}
\begin{proof}
	See \cite{Goetze}, Theorem 2.7.
\end{proof}

One more lemma from \cite{Goetze} describes how close are the distributions of the norms of two Gaussian elements in terms of their covariance operators.
Note that the bound is dimension free.
\begin{lemma}[Gaussian comparison]\label{Lemma: Gaussian comparison}
	Let \( \xi \) and \( \eta \) be Gaussian elements in Hilbert space \( \mathbb{H} \) with zero mean and covariance operators \( \Su_{\xi} \) and \( \Su_{\eta} \), respectively.
	The following inequality holds
\begin{EQA}
	&& \sup_{x \in \R} \left| \Pro( \| \xi \|^2 \geq x ) - \Pro( \| \eta \|^2 \geq x ) \right| \lesssim  \| \Su_{\xi} - \Su_{\eta}\|_1 \times \\
	&& \quad \quad\times
	\left( \frac{1}{\| \Su_{\xi}\|_2^{1/2} \left(\| \Su_{\xi}\|_2^2 - \|\Su_{\xi}\|_{\infty}^2\right)^{1/4}} 
	+ \frac{1}{\| \Su_{\eta}\|_2^{1/2} \left(\| \Su_{\eta}\|_2^2 - \|\Su_{\eta}\|_{\infty}^2\right)^{1/4}}\right) .
\label{Bound: Gaussian comparison}
\end{EQA}
\end{lemma}
\begin{proof}
	See \cite{Goetze}, Theorem 2.1.
\end{proof}

	\section{Auxiliary proofs} \label{Appendix: B}
\subsection{Proof of Theorem \ref{Theorem: clt}}
The proof consists of three steps.

\paragraph{Step 1 }
Apply the representation (\ref{Formula: Projector decomposition}) from Lemma \ref{Lemma: Projector decomposition} to \( \Pe - \Pt \):
\begin{EQA}
    \Pe - \Pt 
    &=& 
    L_{\subspaceset}(\Se-\St) + R_{\subspaceset}(\Se-\St).
\end{EQA}
Then, for \( n\| \Pe - \Pt  \|_{\Fr}^2 \) one has
\begin{EQA} 
    n\| \Pe - \Pt \|_{\Fr}^2 
    & = & 
    n\|L_{\subspaceset}(\Se-\St)\|_{\Fr}^2 
    + n\| R_{\subspaceset}(\Se-\St) \|_{\Fr}^2 \\ 
    && \qquad + 2n \langle L_{\subspaceset}(\Se-\St), R_{\subspaceset}(\Se-\St)\rangle_{2}.
\end{EQA}
Let us estimate how good \( n\|L_{\subspaceset}(\Se-\St)\|_{\Fr}^2 \) approximates \( n\| \Pe - \Pt \|_{\Fr}^2 \): clearly, we have
\begin{EQA}
    && \nquad
    \left| n\| \Pe - \Pt \|_{\Fr}^2 - n\|L_{\subspaceset}(\Se-\St)\|_{\Fr}^2 \right| 
    \\
    & \leq &
    n\| R_{\subspaceset}(\Se-\St)\|_{\Fr}^2 
    + 2n \| L_{\subspaceset}(\Se-\St) \|_{\Fr} \|R_{\subspaceset}(\Se-\St)\|_{\Fr} .
\end{EQA}
Let us elaborate on the right-hand side. First, since
\begin{EQA}
    R_{\subspaceset}(\Se-\St) 
    &=& 
    \Pe - \Pt - \sum_{r \in {\subspaceset}} \sum_{s \notin {\subspaceset}} 
      \frac{\Ptr (\Se-\St) \Pts + \Pts (\Se-\St) \Ptr}{\mus_r - \mus_s}
\end{EQA}
and \( \Pe \), \( \Pt \), \( \sum_{r \in {\subspaceset}}\sum_{s\notin \subspaceset} \Ptr (\Se-\St) \Pts \) have rank at most \( \ms_{\subspaceset} \), then the rank of \( R_{\subspaceset}(\Se-\St) \) is at most \( 4 \ms_{\subspaceset} \). 
Hence, due to the relation between the Frobenius and the operator norms via rank, we have
\begin{EQA}
    \| R_{\subspaceset}(\Se-\St) \|_{\Fr}
    & \leq &
    \sqrt{4\ms_{\subspaceset}}\; \| R_{\subspaceset}(\Se-\St)\|_{\infty}.
\end{EQA}
The bound (\ref{Formula: Bound}) from Lemma \ref{Lemma: Projector decomposition} gives
\begin{EQA}
    \| R_{\subspaceset}(\Se-\St) \|_{\Fr}
    & \leq &
    \sqrt{4\ms_{\subspaceset}} \cdot 15\,\left(1 + \frac{2}{\pi}\frac{l^*_{\subspaceset}}{\gs_{\subspaceset}}\right) \,\frac{\| \Se - \St \|_{\infty}^2}{{\gs_{\subspaceset}}^2}.
\end{EQA}
Now let us bound \( \|L_{\subspaceset}(\Se-\St)\|_{\infty} \): 
\begin{EQA}
    && \| L_{\subspaceset}(\Se-\St) \|_{\infty} 
    =  
    \left\| \sum_{r \in \subspaceset} \sum_{s \notin \subspaceset} \frac{\Ptr(\Se-\St)\Pts + \Pts(\Se-\St)\Ptr}{\mus_r - \mus_s} \right\|_{\infty}  
    \\
    && \qquad \leq
    2 \left\| \sum_{r \in \subspaceset} \sum_{s \notin \subspaceset} \frac{\Ptr(\Se-\St)\Pts }{\mus_r - \mus_s} \right\|_{\infty} 
    =
    2 \left\| \sum_{r \in \subspaceset} \Ptr \sum_{s \notin \subspaceset} \frac{(\Se-\St)\Pts }{\mus_r - \mus_s} \right\|_{\infty}  
     \\
    && \qquad \leq 
     2  \sum_{r \in \subspaceset} \| \Ptr \| \left\| \sum_{s \notin \subspaceset} \frac{\Pts }{\mus_r - \mus_s} \right\|_{\infty} \| \Se-\St\|_{\infty} 
     \\
    && \qquad \leq 
    \frac{2|\subspaceset|\| \Se-\St\|_{\infty}}
    	 {\min\limits_{r\in\subspaceset,\;s\notin\subspaceset} |\mus_r-\mus_s|}  
    \leq 
    2|\subspaceset|\,\frac{\| \Se-\St\|_{\infty}}{\gs_{\subspaceset}} \,.
\end{EQA}
Then, for \( \|L_{\subspaceset}(\Se-\St)\|_{\Fr} \) we have 
\begin{EQA}
    \| L_{\subspaceset}(\Se-\St) \|_{\Fr} &= \sqrt{2\ms_{\subspaceset}} \; \| L_{\subspaceset}(\Se-\St) \|_{\infty}  \leq
     \sqrt{2\ms_{\subspaceset}}\;2|\subspaceset|\,\frac{\| \Se-\St\|_{\infty}}{\gs_{\subspaceset}}\,.
\end{EQA}
Putting this all together, we obtain
\begin{EQA}
    && \nquad \left| n\| \Pe-\Pt \|_{\Fr}^2 - n\|L_r(\Se-\St)\|_{\Fr}^2 \right| 
    \\ 
    & \lesssim &
    n \ms_{\subspaceset} \, \left( 1 + \frac{l^*_{\subspaceset}}{\gs_{\subspaceset}}\right)^2 \,\frac{\|\Se-\St\|_{\infty}^4 }{{\gs_{\subspaceset}}^4}
    + n \ms_{\subspaceset} |\subspaceset| \, \left( 1 + \frac{l^*_{\subspaceset}}{\gs_{\subspaceset}}\right) \,\frac{\|\Se-\St\|_{\infty}^3}{{\gs_{\subspaceset}}^3} .
\end{EQA}
The concentration condition for the sample covariance (\ref{Eq: hat_delta_n}) provides
\begin{equation}
\begin{aligned}
    & \left| n\| \Pe - \Pt \|_{\Fr}^2 - n\|L_{\subspaceset}(\Se-\St)\|_{\Fr}^2 
    \right|  
    \lesssim 
    \overline{\Delta} \,,\\ 
    & \overline{\Delta}
    = 
    n \ms_{\subspaceset} \left( 1 + \frac{l^*_{\subspaceset}}{\gs_{\subspaceset}}\right)\, \left(  
            \left( 1 + \frac{l^*_{\subspaceset}}{\gs_{\subspaceset}}\right)\frac{\hdelta_n^4}{{\gs_{\subspaceset}}^4} \;\lor\; 
            \frac{|\subspaceset| \;\hdelta_n^3}{{\gs_{\subspaceset}}^3} 
    \right) \qquad
\label{Eq: overline_Delta}
\end{aligned}
\end{equation}

with probability \( 1 - \frac{1}{n} \).

\paragraph{Step 2 } 
Following \cite{Naumov_BCSFSPOSC}, we can choose \( \{u_j^{*}\}_{j=1}^p \) as an orthonormal basis in \( \R^p \) and represent 
\( n\| L_{\subspaceset}(\Se-\St) \|_{\Fr}^2 \) as
\begin{EQA}
    && n\| L_{\subspaceset}(\Se-\St) \|_{\Fr}^2 
    =
    n \left\| \sum_{r \in {\subspaceset}} \sum_{s \notin {\subspaceset}} 
            \frac{\Ptr (\Se-\St) \Pts + \Pts (\Se-\St) \Ptr}{\mus_r - \mus_s} \right\|_{\Fr}^2 
    \\
    && \qquad =   
    n \sum_{l,k=1}^p \left({\us_k}^{\T} \;\sum_{r \in {\subspaceset}} 
    \sum_{s \notin {\subspaceset}} 
            \frac{\Ptr (\Se-\St) \Pts + \Pts (\Se-\St) \Ptr}{\mus_r - \mus_s} \;\us_l \right)^2  
    \\
    && \qquad = 
    n \sum_{l,k=1}^p \left({\us_k}^{\T} \,
    \sum_{r_1 \in {\subspaceset}} \,
    \sum_{s_1 \notin {\subspaceset}} 
        \frac{\Projs_{r_1} (\Se-\St) \Projs_{s_1} + \Projs_{s_1} (\Se-\St) \Projs_{r_1}}{\mus_{r_1} - \mus_{s_1}} \;\us_l \right) 
    \\
    && \hspace{1.75cm}
    \times \left({\us_k}^{\T} \;\sum_{r_{2} \in {\subspaceset}} \sum_{s_{2} \notin {\subspaceset}} 
            \frac{\Projs_{r_{2}} (\Se-\St) \Projs_{s_{2}} + \Projs_{s_{2}} (\Se-\St) \Projs_{r_{2}}}{\mus_{r_{2}} - \mus_{s_{2}}} \;\us_l \right) 
    \\
    && \qquad = 
    n \sum_{l,k=1}^p 
            \sum_{\substack{ r_1 \in {\subspaceset} \\ s_1 \notin {\subspaceset}}}
            \sum_{\substack{ r_{2} \in {\subspaceset} \\ s_{2} \notin {\subspaceset}}}
            \left({\us_k}^{\T} \; \frac{\Projs_{r_1} (\Se-\St) \Projs_{s_1} + \Projs_{s_1} (\Se-\St) \Projs_{r_1}}{\mus_{r_1} - \mus_{s_1}} \;\us_l \right)  
    \\
    && \hspace{3.25cm} 
    \times \left({\us_k}^{\T} \; 
            \frac{\Projs_{r_{2}} (\Se-\St) \Projs_{s_{2}} + \Projs_{s_{2}} (\Se-\St) \Projs_{r_{2}}}{\mus_{r_{2}} - \mus_{s_{2}}} \;\us_l \right).
\end{EQA}
As we can see, the  only terms that survive in this sum are the terms with \( r_1 = r_{2} = r \in \subspaceset \), \( s_1 = s_{2} = s \notin \subspaceset \), \( k \in \Deltas_{r} \), \( l \in \Deltas_{s} \), and due to the symmetry the factor \( 2 \) appears. So, we derive 
\begin{EQA}
    n\| L_{\subspaceset}(\Se-\St) \|_{\Fr}^2 
    &=&
    2n \sum_{\substack{k \in \Deltas_r, \\ r \in \subspaceset}} \; \sum_{\substack{l \in \Deltas_s, \\ s \notin \subspaceset}} 
         \left( {\us_k}^{\T} \; \frac{\Ptr (\Se-\St) \Pts}{\mus_{r} - \mus_{s}} \;\us_l \right)^2 
    \\
    &=& 
    2n \sum_{\substack{k \in \Deltas_r, \\ r \in \subspaceset}} \; \sum_{\substack{l \in \Deltas_s, \\ s \notin \subspaceset}} 
         \left( \frac{{\us_k}^{\T} (\Se-\St) \;\us_l}{\mus_{r} - \mus_{s}} \right)^2.
\end{EQA}    
Now let us define for all \( k \in \IndS \) and \( l \notin \IndS \)
\begin{EQA}
    S_{\subspaceset}(u_k^{*}, u_l^{*}) 
    &=& 
    \sqrt{2n} \; \frac{{\us_k}^{\T} (\Se-\St) \;\us_l}{\mus_{r} - \mus_{s}}.
\end{EQA}
This set of quantities can be considered as matrix
\begin{EQA}[c]
    \{ S_{\subspaceset}(u_k^{*}, u_l^{*}) \}_{\substack{k \in \IndS \\ l \notin \IndS}} \in \R^{\ms_{\subspaceset} \times (p-\ms_{\subspaceset})},
\end{EQA}
or, we can arrange a vector \( S_{\subspaceset} \in \R^{\ms_{\subspaceset}(p-\ms_{\subspaceset})} \) with components \( S_{\subspaceset}(u_k^{*}, u_l^{*}) \) ordered in some particular way.
Let us notice that
\begin{EQA}
    n\| L_{\subspaceset}(\Se-\St) \|_{\Fr}^2 
    &=& 
    \| S_{\subspaceset} \|^2.
\end{EQA}

\paragraph{Step 3 }
Now our goal is to show that  \( S_{\subspaceset} \) is approximately \( \ND(0, \GammasJ) \) using a version of Berry-Esseen theorem given by \cite{Bentkus_LTB}. Represent \( S_{\subspaceset} \) as
\begin{EQA}
     S_{\subspaceset} 
     &=& 
     \frac{1}{\sqrt{n}} \sum_{j=1}^n S^{(j)}, 
\end{EQA}
where \( S^{(j)} \) is a random vector with components
\begin{EQA}
    S^{(j)}(u_k^{*}, u_l^{*}) 
    &=& 
    \frac{\sqrt{2}}{\mus_{r} - \mus_{s}}
    	\;({\us_k}^{\T} X_j) \cdot ({\us_l}^{\T} X_j)
\end{EQA}
for all \( k \in \IndS \) and \( l \notin \IndS \). 

It is straightforward to verify that the covariance matrix of \( S^{(j)} \) (and hence of \( S_{\subspaceset} \)) is \( \GammasJ \) from (\ref{Formula: Gamma_r}) under the condition that \( \Pt X_j \) and \( (\I_{p} - \Pt)X_j \) are independent. 
Consider an entry of the covariance matrix of \( S^{(j)} \) indexed by \( (k, l) \) and \( (k^{\prime}, l^{\prime}) \), 
where \( k \in \Deltas_r, k^{\prime} \in \Deltas_{r^{\prime}}, \;r, r^{\prime} \in \subspaceset \) 
and \( l \in \Deltas_s, l^{\prime} \in \Deltas_{s^{\prime}}, \;s, s^{\prime} \notin \subspaceset \):
\begin{EQA}
    Cov\left(S^{(j)}\right)_{\substack{(k,l) \\ (k^{\prime}, l^{\prime})}} 
    &=& 
    \E \left( S^{(j)}(u_k^{*}, u_l^{*}) \cdot S^{(j)}(u_{k^{\prime}}^{*}, u_{l^{\prime}}^{*}) \right) 
    \\
    &=& 
    \frac{2\;\E \left[({\us_k}^{\T} X_j) \cdot ({\us_l}^{\T} X_j) \cdot ({\us_{k^{\prime}}}^{\T} X_j) \cdot ({\us_{l^{\prime}}}^{\T} X_j)\right]}{(\mus_{r} - \mus_{s})(\mus_{r^{\prime}} - \mus_{s^{\prime}})}.
\end{EQA}
Now, the independence of \( \Pt X_j \) and \( (\I_{p} - \Pt)X_j \) implies the independence of \( (u_k^{*}, u_{k^{\prime}}^{*})^{\T} \Pt X_j \) and \( (u_l^{*}, u_{l^{\prime}}^{*})^{\T} (\I_{p} - \Pt)X_j \), which can be rewritten as independence of \( ({u_k^{*}}^{\T} X_j, {u_{k^{\prime}}^{*}}^{\T} X_j)^{\T} \) and \( ({u_l^{*}}^{\T} X_j, {u_{l^{\prime}}^{*}}^{\T} X_j)^{\T} \).
This means that the expectation in the expression for the covariance entry can be splitted as
\begin{EQA}
    Cov\left(S^{(j)}\right)_{\substack{(k,l) \\ (k^{\prime}, l^{\prime})}} 
    &=& 
    \frac{2\; \E \left[({\us_k}^{\T} X_j)({\us_{k^{\prime}}}^{\T} X_j)\right] 
    	  \cdot \E\left[({\us_l}^{\T} X_j) ({\us_{l^{\prime}}}^{\T} X_j)\right]}
    	 {(\mus_{r} - \mus_{s})(\mus_{r^{\prime}} - \mus_{s^{\prime}})}.
\end{EQA}
The observation  that \( {\us_k}^{\T} \St \us_{k^{\prime}} = \mus_r \cdot \Ind\{ k=k^{\prime}\} \) and \( {\us_l}^{\T} \St \us_{l^{\prime}} = \mus_s \cdot \Ind\{ l=l^{\prime}\} \) establishes the fact that \( Cov(S^{(j)}) = \GammasJ \).

To apply Theorem 1.1 from \cite{Bentkus_LTB}, we need to bound \( \E \| {\GammasJ}^{-1/2} S^{(j)}\|^3 \).
First, let us notice that 
\begin{EQA}
    \left[{\GammasJ}^{-1/2} S^{(j)}\right](u_k^{*}, u_l^{*}) 
    &=& 
    \frac{{\us_k}^{\T} X_j}{\sqrt{\mus_r}} \cdot 
    \frac{{\us_l}^{\T} X_j}{\sqrt{\mus_s}}.
\end{EQA}
Further, recalling the auxiliary matrices
\begin{EQ}[rcl]
    \Ut & \eqdef & 
    \left\{ \frac{{\us_{k}}^{\T}}{\sqrt{\mus_{r}}} \right\}_{k\in\IndS } 
    \in \R^{\ms_{\subspaceset} \times p},
    \\
    \Vt 
    & \eqdef &
    \left\{ \frac{{\us_{l}}^{\T}}{\sqrt{\mus_{s}}} \right\}_{l\notin\IndS} 
    \in \R^{(p- \ms_{\subspaceset})  \times p} \,,
\label{Equation: UV}
\end{EQ} 
we have
\begin{EQA}
    &&\nquad \| {\GammasJ}^{-1/2} S^{(j)} \|^2 = 
    \sum_{k\in\IndS} \sum_{l\notin\IndS}
    \frac{({\us_k}^{\T} X_j)^2}{\mus_r} \cdot \frac{({\us_l}^{\T} X_j)^2}{\mus_s} 
    \\
  &=&
  \left\{ \sum_{k\in\IndS} \frac{({\us_k}^{\T} X_j)^2}{\mus_r}  \right\} \cdot
  \left\{ \sum_{l\notin\IndS} \frac{({\us_l}^{\T} X_j)^2}{\mus_s}  \right\} 
  =
  \| \Ut X_j \|^2 \;\| \Vt X_j \|^2.
\end{EQA}
Then, 
\begin{EQA}
    && \nquad \E \| {\GammasJ}^{-1/2} S^{(j)} \|^3
    = \E \left( \| \Ut X_j \|^3 \;\| \Vt X_j \|^3 \right) 
    = \E  \| \Ut X_j \|^3 \cdot \E \| \Vt X_j \|^3,
\end{EQA}
where we again used the fact that the independence of \( \Pt X_j \) and \( (\I_{p} - \Pt)X_j \) implies the independence of \( \Ut \Pt X_j = \Ut X_j \) and \( \Vt (\I_{p} - \Pt)X_j = \Vt X_j \).

Therefore, Theorem 1.1 from \cite{Bentkus_LTB} yields
\begin{EQA}
    \sup_{x \in \R} 
    \left| \Pro\left( \| S_{\subspaceset}\|^2 \leq x\right) 
    	- \Pro(\| \xi\|^2 \leq x) 
    \right|
    & \lesssim &
    \E  \| \Ut X  \|^3 \cdot \E \| \Vt X \|^3 \cdot \frac{p^{1/4}}{\sqrt{n}},
\end{EQA}
or, recalling that \( \| S_{\subspaceset}\|^2 = n\| L_{\subspaceset}(\Se - \St) \|^2_{2} \),
\begin{EQA}
    && \nquad \sup_{x \in \R} \left| 
    	\Pro\left( n\| L_{\subspaceset}(\Se - \St) \|^2_{2} \leq x\right) 
		- \Pro(\| \xi\|^2 \leq x) 
    \right|\\
    & \lesssim &
    \E  \| \Ut X \|^3 \cdot \E \| \Vt X \|^3 \cdot \frac{p^{1/4}}{\sqrt{n}} \, .
\end{EQA}

\paragraph{Step 4 }
Next, for \( \overline{\Delta} \) defined by (\ref{Eq: overline_Delta}) from Step \( 1 \) we may write for any \( x\in\R \)
\begin{EQA}
      \Pro\left( n\| \Pe - \Pt \|_{\Fr}^2 \geq x\right) 
      & \leq &
      \Pro\left( n\| L_{\subspaceset}(\Se-\St)\|_{\Fr}^2 \geq x - \overline{\Delta}\right) 
      \\
      && \quad 
      + \; \Pro\left( 
      	n\| \Pe - \Pt \|_{\Fr}^2 - n\| L_{\subspaceset}(\Se-\St)\|_{\Fr}^2 
      		\geq \overline{\Delta} 
	\right).
\end{EQA}
Hence, 
\begin{EQA}
      &&\nquad 
      \sup_{x\in\R} \left\{\Pro\left( n\| \Pe-\Pt \|_{\Fr}^2 \geq x \right) - \Pro\left( \|\xi\|^2 \geq x\right) \right\} 
      \\
      & \leq &
      \sup_{x\in\R} \left\{ \Pro\left( \| L_{\subspaceset}(\Se-\St)\|^2 \geq x -\overline{\Delta}\right)
                                            - \Pro\left( \|\xi\|^2 \geq x-\overline{\Delta}\right) \right\} 
      \\
      && \quad + \; \sup_{x\in\R} \left\{ \Pro\left( \|\xi\|^2 \geq x-\overline{\Delta}\right)
                                            - \Pro\left( \|\xi\|^2 \geq x\right) \right\} 
      \\
      && \quad + \; \Pro\left( n\| \Pe-\Pt \|_{\Fr}^2 - n\| L_{\subspaceset}(\Se-\St)\|_{\Fr}^2 \geq \overline{\Delta} \right).
\end{EQA}
The first term in the right-hand side was bounded in Step \( 3 \) by 
\[ 
	\E  \| \Ut X \|^3 \cdot \E \| \Vt X \|^3 \cdot  \frac{p^{1/4}}{\sqrt{n}} .
\]
The second term is bounded by \( \frac{\overline{\Delta}}{\| \Gammas_{\subspaceset} \|_2^{1/2} \left( \| \Gammas_{\subspaceset} \|_2^2 - \| \Gammas_{\subspaceset} \|_{\infty}^2 \right)^{1/4}} \) according to the Anti-concentration Lemma~\ref{Lemma: Anticoncentration}. 
The last term is less than \( {1}/{n} \) in view of (\ref{Eq: overline_Delta}) from Step 1.
Therefore,
\begin{EQA}
      && \nquad
      \sup_{x\in\R} \left\{\Pro\left( n\| \Pe-\Pt \|_{\Fr}^2 \geq x\right) - \Pro\left( \|\xi\|^2 \geq x \right) \right\} 
      \\
      &\lesssim &
      \E  \| \Ut X \|^3 \cdot \E \| \Vt X \|^3 \cdot \frac{p^{1/4}}{\sqrt{n}} + \frac{\overline{\Delta}}{\| \Gammas_{\subspaceset} \|_2^{1/2} \left( \| \Gammas_{\subspaceset} \|_2^2 - \| \Gammas_{\subspaceset} \|_{\infty}^2 \right)^{1/4}} + \frac{1}{n}.
\end{EQA}
Similarly, one can verify that
\begin{EQA}
      &&\nquad 
      \sup_{x\in\R} \left\{ \Pro\left( \|\xi\|^2 \geq x\right) - \Pro\left( n\| \Pe-\Pt \|_{\Fr}^2 \geq x\right) \right\} 
      \\
      & \lesssim &
      \E  \| \Ut X \|^3 \cdot \E \| \Vt X \|^3 \cdot \frac{p^{1/4}}{\sqrt{n}} 
      + \frac{\overline{\Delta}}{\| \Gammas_{\subspaceset} \|_2^{1/2} \left( \| \Gammas_{\subspaceset} \|_2^2 - \| \Gammas_{\subspaceset} \|_{\infty}^2 \right)^{1/4}} 
      + \frac{1}{n}.
\end{EQA}
Putting together the previous two bounds, we derive the final result:
\begin{EQA}
      && \nquad
      \sup_{x\in\R} \left| \Pro\left( n\| \Pe-\Pt \|_{\Fr}^2 \geq x\right) - \Pro\left( \|\xi\|^2 \geq x\right) \right| 
      \\
      & \lesssim &
      \E  \| \Ut X \|^3 \cdot \E \| \Vt X \|^3 \cdot \frac{p^{1/4}}{\sqrt{n}} 
      + \frac{\overline{\Delta}}{\| \Gammas_{\subspaceset} \|_2^{1/2} \left( \| \Gammas_{\subspaceset} \|_2^2 - \| \Gammas_{\subspaceset} \|_{\infty}^2 \right)^{1/4}}
      + \frac{1}{n}.
\end{EQA}